\documentclass[10pt,psamsfonts]{amsart}
\usepackage{amsfonts}
\usepackage{amssymb}
\usepackage{amstext}

\newcommand{\R}{\mathbb R}

\newcommand{\N}{\mathbb N}

\newcommand{\del}{\partial}
\newcommand{\e}{\varepsilon}
\newcommand{\ptl}{\partial}

\newcommand{\Om}{\Omega}

\renewcommand{\leq}{\leqslant}
\renewcommand{\geq}{\geqslant}
\newcommand{\cal}{\mathcal}

\def\ds{\displaystyle}

\newtheorem{theorem}{Theorem}[section]
\newtheorem{lemma}[theorem]{Lemma}
\newtheorem{proposition}[theorem]{Proposition}
\newtheorem{corollary}[theorem]{Corollary}
\newtheorem{definition}[theorem]{Definition}
\newtheorem{remark}[theorem]{Remark}

\title[CMC and the gradient theory of phase transitions]{From constant
mean curvature hypersurfaces \\ to the gradient theory of phase
transitions}

\author[F.~Pacard]{Frank Pacard}
\address{Centre de Math\'ematiques \\ Facult\'e de Sciences et
Technologie \\ Universit\'e Paris XII - Val de Marne \\ 94 010
CRETEIL Cedex \\ France}
\email{pacard@univ-paris12.fr}
\author[M.~Ritor\'e]{Manuel Ritor\'e}
\address{Departamento de Geometr\'{\i}a y Topolog\'{\i}a \\ Facultad
de Ciencias \\ Universidad de Gra\-nada \\ E--18071 Granada \\ Espa\~na}
\email{ritore@ugr.es}

\begin{document}

\maketitle

\thispagestyle{empty}

\section{Introduction}

Let $\Om \subset {\R}^{n+1}$, $n\geq 1$, be an open bounded set with
smooth boundary $\del\Om$.  For any $\e >0$ and any function $u :
\Om\to\R$ such that $u\in H^1(\Om)$, we consider the energy
\[
E_\e (u) := \e^2 \, \int_\Om |\nabla u|^2 \, dx + \int_{\Om}(1-u^2)^2\, dx ,
\]
being understood that $E_\e (u)= \infty$ if $u \notin L^4(\Om)$. We
also consider the constraint
\[
V(u) : =  \int_\Om u \, dx .
\]

Given $c_0 \in (- 1,1)$, we are interested in the critical points of
$E_\e$ subject to the constraint $V(u) = c_0 \, |\Omega|$, where
$|\Omega|$ denotes the volume of $\Omega$.  Any critical point of
this variational problem is solution of
\begin{equation}
\left\{
\begin{array}{rlllllll}
\ds \e^2 \, \Delta u + 2 \, (u -u^3)  & =  &   \e\,\lambda, \quad
&\mbox{in} \quad \Om \\[3mm]
\ds\del_{\nu} u & = & 0,   \quad &\mbox{on}\quad \del \Om ,
\end{array}
\right.
\label{eq:**}
\end{equation}
where $\nu$ denotes a unit vector field normal to $\del\Om$ and where
$\e\,\lambda \in {\R}$ corresponds to the Lagrange multiplier
associated to the constraint $V(u)=c_0 \, |\Omega|$. One can also
ignore the volume constraint, in which case a critical point would
satisfy equation \eqref{eq:**} with $\lambda=0$.

Since classical methods of the calculus of variation apply, there is
no difficulty in finding {\it minimizers} of $E_\e$.  The real issue
is the study of the asymptotic behavior of the minimizers (or more
generally of the critical points) of $E_\e$ as the parameter $\e$
tends to $0$. There has been a number of important work on this
question over the last two decades and the basic result can be
described as follows~: Assume that $(\e_k)_{k \geq 0}$ tends to $0$
and let $(u_k)_{k\geq 0}$ be a sequence of {\em minimizers\/} of
$E_{\e_k}$ under the constraint $V(u)=c_0 \, |\Omega|$.  Then, up to
a subsequence, one can assume that $(|u_k|)_{k\geq 0}$ converges a.e.
to the constant function $1$. In the definition of the energy
$E_{\e}$, the role of the term
\[
\int_\Om (1-u^2)^2 \, dx ,
\]
is precisely to force the sequence of function $(|u_k|)_{k\geq 0}$ to
converge to $1$ when the parameter $\e_k$ tends to $0$. Extracting
subsequences if this is necessary, we can define $\Om^+$
(resp.~$\Om^-$) to be the set of points where $u_k$ converges to $+1$
(resp.~$-1$).  The subsets $\Om^\pm$ are not arbitrary since the
constraint $V(u_k) =c_0 \, |\Omega |$ forces $\Om^\pm \subset \Om$ to
satisfy
\[
|\Om^+|-|\Om^-|= c_0 \, |\Omega| .
\]
Now, the role of the Dirichlet integral
\[
\e^2 \, \int_\Om |\nabla u|^2 \, dx ,
\]
in the definition of $E_\e$ forces the interface between the subsets
$\Om^+$ and $\Om^-$ to be ``as small as possible'', since this is
where the gradient of the function $u_k$ will concentrate when $\e_k$
tends to $0$.  More precisely
\[
N : = \del \Om^+ \cap \Om = \del \Om^- \cap \Om ,
\]
can be shown to be a minimizer of the isoperimetric problem: Minimize
amongst all domains $D \subset \Om$ the $n$-dimensional Hausdorff
measure $\mathcal{H}^n(\del D)$ of the boundary $\del D$ subject to
the volume constraint
\[
|D| = \frac{ 1+ c_0 }{2} \, |\Om| .
\]
We refer to \cite{modica1}, \cite{modica2}, \cite{luckmodica},
\cite{anzbalorl}, \cite{baldo}, \cite{sternberg} for more precise
statements. From a purely analytic point of view, $N$ can be
understood as the limit of the nodal sets of the functions $u_k$, as
$k$ tends to $+\infty$.

\section{Statement of the problem}\setcounter{equation}{0}

It is interesting to generalize the above problem first by
considering instead of $\Om \subset {\R}^{n+1}$, any compact
Riemannian manifold with or without smooth boundary and also by
replacing the nonlinearity $(1-u^2)^2$ by a more general one.

Hence, in this paper, we consider $(M, g)$ to be a
$(n+1)$-dimensional compact Riemannian manifold with or without
smooth  boundary. In the case where $\del M$, the boundary of $M$, is
not  empty, we can assume without loss of generality that $M$ is a
subdomain of a larger Riemannian manifold $(\tilde M, \tilde g)$,
with  $\tilde g_{|M} = g$.  In particular, $\del M$ is a smooth
hypersurface of $\tilde M$.

Let $W : {\mathbb R}\to {\mathbb R}$ be a smooth function which is
positive away from $u= \pm 1$.  We assume~that
\begin{equation}
W(\pm 1) =0,
\label{eq:2.1}
\end{equation}
so that the infimum of $W$ is achieved at the points $u=\pm 1$.
Further assume that these points are nondegenerate critical points of
$W$.  In other words
\begin{equation}
W''(\pm 1) > 0 .
\label{eq:2.2}
\end{equation}

For any $\e >0$ and any function $u : M \longrightarrow {\R}$, such
that $u \in H^1(M)$, we define the energy
\begin{equation}
E_\e (u) := \e^2 \, \int_M |\nabla u|^2_g \, dv_g + \int_{M} W(u) \, dv_g ,
\label{eq:energy}
\end{equation}
where $\nabla$ denotes the gradient and $dv_g$ the volume form on $M$
associated to the Riemannian metric $g$.  As usual, we agree that
$E_\e (u) = \infty$ when $W(u) \notin L^1(M)$. We also define the
volume constraint
\begin{equation}
V(u) : =  \int_M u \, dv .
\label{eq:constraint}
\end{equation}

Granted the above definition, there are two closely related
variational problems we can consider.

\noindent
1 - We consider the critical points of the energy $u \longrightarrow
E_\e(u)$, which are solutions~of
\begin{equation}
\ds - \e^2 \, \Delta_g u + \frac{1}{2} \, W'(u)  = 0,
\label{eq:2.113}
\end{equation}
in $M$, where $\Delta_g$ is the Laplace-Beltrami operator in $M$.
Moreover, if $\del M\neq\emptyset$ then the additional condition
\begin{equation}
       \del_{\nu_{\del M}} u = 0 ,
       \label{eq:2.3333}
\end{equation}
must hold on $\del M$, where $\nu_{\del M}$ denotes the unit vector
field normal to $\del M$. This problem is related to the Allen-Cahn
equation \cite{All-Cah} and it is well known that, as $\e$ tends to
$0$, the interfaces (i.e. the nodal sets of the solutions of
\eqref{eq:2.113}) converge to minimal hypersurfaces. Concerning this
variational problem, the question we would like to address in this
paper is the following~:

\medskip

(P-1) \hspace{10mm} \parbox{100mm}{\em Assume that $N \subset M$ is a
minimal hypersurface. Does $N$ appear as the limit, as the parameter
$\e$ tends to $0$, of the nodal sets of a sequence of critical points
of $E_\e$ ?}

\medskip

\noindent
2 - Given $c_0 \in (- 1,1)$, we consider the critical points of the
energy $u \longrightarrow E_\e(u)$ under the constraint $V (u) = c_0
\, |M|$, where $|M|$ denotes the volume of $M$. This time, such a
critical point $u$ is a solution of
\begin{equation}
\ds - \e^2 \, \Delta_g u + \frac{1}{2} \, W'(u)  =  \e \, \lambda ,
\label{eq:2.3}
\end{equation}
in $M$, where $\e \, \lambda \in {\R}$ corresponds to the Lagrange
multiplier associated to the constraint $V(u)=c_0 \, |M|$. Moreover,
$u$ satisfies (\ref{eq:2.3333}) on $\del M$ if $\del M\neq\emptyset$.
According to \cite{vdw}, \cite{ch}, the energy $E_{\e}$ corresponds
to the total energy of a fluid within the Wan der Waals-Cahn-Hilliard
theory of phase transitions.  The Lagrange multiplier $\e \,
\lambda$, which appears in (\ref{eq:2.3}), is known in the physics
literature as the chemical potential of the density configuration
$u$. Now, the question we would like to address becomes~:

\medskip

(P-2) \hspace{10mm} \parbox{100mm}{\em  Assume that $N \subset M$ is
a constant mean curvature hypersurface.  Does $N$ appear as the
limit, as the parameter $\e$ tends to $0$, of the nodal sets of a
sequence of critical points of $E_\e$ subject to the constraint $V(u)
= c_0 \, |M|$ ?}

\medskip

Before we proceed, let us observe that, in both problems, we are not
only looking for minimizers of $E_\e$ but more generally for critical
points.

\begin{remark}
If the infimum of the function $W$ is achieved at exactly two points
$u_\pm$, there is no loss of generality in considering that $u_\pm =
\pm 1$ since we can always reduce to this case by considering $u
\mapsto W(au + b)$ where $a$ and $b$ are chosen appropriately.
\end{remark}

\section{Definitions and Preliminaries}\setcounter{equation}{0}

\subsection{Admissible hypersurfaces in $M$}

Obviously if $N$ is the nodal set of some function $u$ which is
defined in $M$ and if $0$ is a regular value of $u$ then $M - N$ is
the union of
\begin{equation}
M^+(N) : = u^{-1}((0, + \infty)) \qquad \mbox{and} \qquad M^{-} (N) :
= u^{-1}((-\infty, 0)).
\label{eq:3.1}
\end{equation}
We shall associate to $N$ the unit normal vector field which points
into $M^+(N)$. In the case where $M$ has a boundary, it may happen
that $N$ also has a boundary $\del N \subset \del M$.  In this case,
if $N$ is the nodal set of the function $u$ and if in addition the
function $u$ has $0$ Neumann boundary condition on $\del M$, then for
all $p \in \del N \subset \del M$, the normal vector to $N$ at $p$
and the normal vector to $\del M$ at $p$ are orthogonal. This later
condition is standard in the study of minimal and constant mean
curvature hypersurfaces. Indeed, it is well known that smooth
hypersurfaces $N$ which are stationary points of the area functional
(possibly with a volume constraint) and have a boundary $\del
N\subset\del M$, satisfy the later orthogonality condition. This
motivates the following~:

\begin{definition} A smooth embedded hypersurface $N \subset M$ (not
necessarily connected) is {\em admissible\/} if $N$ is the nodal set
of a smooth function $u$ for which $0$ is a regular value of $u$ and
which, in the case where $M$ has a boundary, has $0$ Neumann boundary
condition.
\label{de:3.1}
\end{definition}

A hypersurface $N\subset M$ which separates $M$ into two regions
$M^\pm (N)$, and which meets $\del M$ orthogonally in the case where
$N$ has a nonempty boundary, is easily shown to be admissible by
using partitions of unity.

\subsection{The Jacobi operator}

Before we introduce our next definition, we recall a few basic facts
about the study of constant mean curvature hypersurfaces $N$ in a
Riemannian manifold $(M,g)$. To begin with, let us recall that the
Jacobi operator, that is the linearized mean curvature operator about
$N$, is given by
\begin{equation}
{\cal L}_N : = \Delta_N + |A_N|^2 + \mbox{Ric}_g (\nu_N, \nu_N),
\label{eq:3.3}
\end{equation}
where $\Delta_N$ is the Laplace-Beltrami operator on $N$, $|A_N|^2$
denotes the norm of the second fundamental form of $N$,
$\mbox{Ric}_g$ is the Ricci tensor of $M$ and $\nu_N$ is a unit
normal to $N$.

Given any (smooth) small function $w$ on $N$, we can consider the
hypersurface $N(w)$, the normal graph on $N$ of the function $w$ (the
image of $N$ by the map $p\in N\mapsto \exp_p(w(p)\nu_N(p))$). If
$H(w)$ denotes the mean curvature of $N (w)$, defined as the
arithmetic mean of the principal curvatures, then the linear operator
${\cal L}_N$ is the differential of $w \mapsto nH(w)$ at $w \equiv 0$.

When $\partial N$ is empty, solutions of the homogeneous problem
\begin{equation}
{\cal L}_N\, w = 0,
\end{equation}
on $N$ are called Jacobi fields.  When $\partial N$ is not empty, we
further assume that $N$ meets $\del M$ orthogonally, then Jacobi
fields are the solutions of ${\cal L}_N \, w =0$ in $N$ which satisfy
the boundary condition
\begin{equation}
{\cal B}_N  \, w : =  \partial_{\nu_{\del M}} \, w + A_{\del M}
(\nu_N , \nu_N) \, w = 0,
\label{eq:3.4}
\end{equation}
on $\del N$, where $A_{\del M}$ is the second fundamental form of
$\partial M$ in $\tilde M$.  Equation (\ref{eq:3.4}) has its origin
in the requirement that all the hypersurfaces we are looking at meet
$\del M$ orthogonally and this should be true for the hypersurfaces
generated by the flow associated to $X$.

Minimal hypersurfaces are critical points of the area functional
while constant mean curvature hypersurfaces are critical points of
the area functional with respect to deformations that keep constant
the volume enclosed by the hypersurface. Consider a deformation of
the hypersurface $N$ by the flow generated by a vector field $X$. The
second variation formula for the area functional is then given by
\[
X \longrightarrow - \, \int_N w \, \mathcal{L}_N w \, d a_g +
\int_{\partial N} w \, {\cal B}_N \, w \, ds_g,
\]
where the function $w : = g(\nu_N, X)$. If one considers a
deformation which is volume preserving up to first order, then the
function $w$ also has to satisfy
\[
\int_N w \, d a_g = 0.
\]
We refer to \cite{bdce} or \cite{rossouam} for a derivation of the
second variation of the area functional in a Riemannian manifold.
Here $da_g$ and $ds_g$ are the volume forms on $N$ and $\del N$ which
are induced by the metric $g$.

\subsection{Nondegeneracy}

The previous definitions being understood, we can now give the
notions of nondegeneracy which are associated to the two problems we
are interested in. To begin with let us define the notion of
nondegenerate minimal hypersurface~:
\begin{definition}
An admissible minimal hypersurface $N$ is said to be nondegenerate if
there are no nontrivial solutions $w \in {\cal C}^{2, \alpha}(N)$ of
\[
\mathcal{L}_N \, w =0,
\]
in $N$, with ${\cal B}_N \, w =0$ on $\del N$ if $N$ has a boundary.
\label{de:3.4}
\end{definition}
The notion of nondegeneracy for minimal hypersurfaces is standard.
Consider the Jacobi operator
\begin{equation}
{\cal L}_N :  \left[ {\cal C}^{2, \alpha}(N) \right]_0
\longrightarrow  {\cal C}^{0, \alpha}(N),
\label{eq:wwww}
\end{equation}
where the subscript $0$ is meant to point out that functions in
$\left[ {\cal C}^{2, \alpha}(N)\right]_0$ satisfy ${\cal B}_N w=0$ on
$\del N$ when this latter is not empty.  Nondegeneracy is equivalent
to the fact that the operator ${\cal L}_N$ is injective. This
operator being self-adjoint and elliptic, nondegeneracy is also
equivalent to the invertibility of the operator ${\cal L}_N$ defined
in (\ref{eq:wwww}). On a more geometric point of view, if $N$ is a
nondegenerate minimal hypersurface, the implicit function theorem
ensures that it is possible to find a hypersurface $\tilde N$ which
is close to $N$ and whose mean curvature $\tilde H$ is prescribed,
close to the mean curvature $H$ of $N$.

We will also need the notion of volume-nondegenerate constant mean
curvature hypersurface~:
\begin{definition}
An admissible constant mean curvature hypersurface $N$ is said to be
volume-nondegenerate if there are no nontrivial solutions $(w,c) \in
{\cal C}^{2, \alpha}(N) \times {\mathbb R}$ of
\[
\mathcal{L}_N \, w  + c =0, \qquad  \mbox{and} \qquad  \int_N w \, dv_g = 0,
\]
in $N$, with ${\cal B}_N \, w =0$ on $\del N$ if $N$ has a boundary.
\label{de:3.44}
\end{definition}
The notion of volume-nondegeneracy is less standard and perhaps
requires some explanation. This time, we consider the
extended-operator
\begin{equation}
\begin{array}{rlrlll}
L_N & : & \left[ {\cal C}^{2, \alpha}(N) \right]_0 \times {\mathbb R}
& \longrightarrow & \quad {\cal C}^{0, \alpha}(N) \times {\mathbb
R}\\[3mm]
& & (w, c)  & \longmapsto & \displaystyle \left( {\cal L}_N w + c,
\int_N w  \, da_g \right).
\end{array}
\label{eq:extop}
\end{equation}
This definition being understood, volume-nondegeneracy is equivalent
to the fact that the  operator $L_N$ is injective. Observe that $L_N$
is self-adjoint with respect to the scalar product
\[
\langle (v,c), (w, d) \rangle : =  \int_N v\, w \, da_g + c \, d ,
\]
in $L^2(N)\times {\mathbb R}$.  The operator $L_N$ being clearly
elliptic, volume-nondegeneracy is also equivalent to the
invertibility of the operator $L_N$ defined in (\ref{eq:extop}).
>From a geometric point of view, if $N$ is a constant mean curvature
volume-nondegenerate hypersurface, the implicit function theorem
ensures that it is possible to find a hypersurface $\tilde N$ which
is close to $N$ and whose mean curvature $\tilde H$ is, {\em up to a
constant function}, prescribed close to $H$ the mean curvature of $N$
and such that the volume enclosed by this hypersurface $M^+(\tilde
N)$ is prescribed close to $M^+ (N)$, the volume enclosed by $N$.
Hence, it is possible to prescribe the volume enclosed by $\tilde N$
and, up to a constant function, the mean curvature of $\tilde N$.

\section{Statement of the result}\setcounter{equation}{0}

The previous definitions being understood, we can now state the
results we have obtained concerning both (P-1) and (P-2).

We have the~:
\begin{theorem}
Assume that $N \subset M$ is an admissible nondegenerate minimal
hypersurface. Then, there exists $\e_0 >0$ and for all $\e \in (0,
\e_0)$ there exists $u_\e$, critical point of $u \longrightarrow E_\e
(u)$, such that $u_\e$ converges uniformly to $1$ on compacts subsets
of $M^+ (N)$ (resp.  to $-1$ on compacts subsets of $M^{-} (N)$).
\label{th:4.11}
\end{theorem}
Let us mention the work of M. Kowalczyk \cite{Kow} where a similar result is obtained when $M$ is a two dimensional domain of ${\R}^2$ and $N$ is a line segment. 

We will also prove the~:
\begin{theorem}
Assume that $N \subset M$ is an admissible volume-nondegenerate
constant mean curvature hypersurface. Then, there exists $\e_0 >0$
and for all $\e \in (0, \e_0)$ there exists $u_\e$, critical point of
$u \longrightarrow E_\e (u)$ under the constraint $V (u) = |M^+
(N)|-|M^-(N)|$, such that $u_\e$ converges uniformly to $1$ on
compacts subsets of $M^+ (N)$ (resp.  to $-1$ on compacts subsets of
$M^{-} (N)$).
\label{th:4.1}
\end{theorem}

It is in general extremely hard to check whether a given minimal
hypersurface (resp. constant mean curvature hypersurface) is
nondegenerate (resp. volume-nondegenerate). Hopefully, first observe
that both nondegeneracy and volume-nondegeneracy are ``open
conditions", namely are stable under small perturbation of the
metric. Moreover, in \cite{white}, B. White has proved that minimal
hypersurfaces are nondegenerate for a generic choice of the metric.
It follows from similar arguments that volume-nondegeneracy also
holds for a generic choice of the metric and of the mean curvature.

The solutions constructed in Theorem~2 are solutions of
\[
-\e^2\Delta u_{\e}+\frac{1}{2} \,W'(u_{\e})= \e \, \lambda_{\e}.
\]
As a byproduct of our construction, we obtain a precise expansion of
$u_\e$ in terms of $\e$. In particular, we get the expansion of the
Lagrange multiplier $\lambda_\e$
\[
\lambda_{\e} = \frac{1}{2} \, c_\star \,n \, H_N + {\cal O}(\e),
\]
where $H_N$ is the mean curvature of the limit interface $N$ and
where the constant $c_\star$ is given by
\[
c_\star : = \int_{-1}^{+1} \sqrt {W (s)} \,ds.
\]
Finally, in both problems, the expansion of the energy $E_\e (u_\e)$
of the solutions we construct, is given by
\[
E_\e (u_\e) = 2 \, \e \, c_\star  \, |N| + {\cal O} (\e^2).
\]
where $|N|$ is the volume of the interface $N$. These expansions
agree with the expansions which have already been obtained in
\cite{luckmod} in the case where $u_\e$ are {\em minimizers} of
$E_\e$ subject to the constraint $V =  |M^+ (N)|-|M^-(N)|$.

Unfortunately, in many interesting cases and despite the genericity
of these notions, minimal hypersurfaces are degenerate and constant
mean curvature hypersurfaces are volume-degenerate. This is for
example the case when there is a nontrivial group of isometries
acting on $M$. It is well known that any $(\phi_t)_{t\in (-1, 1)}$
smooth one-parameter group of isometries of $M$ which contains the
identity, for example $\phi_0 = \text{Id}$, gives rise to a Jacobi
field on $N$ (when $M$ has a non empty boundary, we ask that these
isometries preserve globally $\partial M$). Actually, the Jacobi
field $w$ is explicitly given by $w := g(\nu_N ,X)$, where $\nu_N$ is
the normal vector field to $N$ and where $X : = {\del_t
\phi_t}_{|t=0}$ is the Killing field corresponding to the
one-parameter group of isometries $\{\phi_t\}_{t\in (-1, 1)}$.
Observe that the isometries $\phi_t$ preserve the volume of the
regions $M^\pm (N)$. Therefore, it follows from the first variation
of volume that the Jacobi field $w$ has mean zero on $N$. In
particular, $w$ is a nontrivial solution of ${\cal L}_N w = 0$ (resp.
$(w,0)$ is a nontrivial solution of $L_N (w,0)=0$) and the
hypersurface $N$ is degenerate (resp. volume-degenerate).

In some cases, it is possible to reduce to a nondegenerate (or
volume-nondegene\-ra\-te) problem by working in the space of
functions and hypersurfaces which are equivariant with respect to
some finite group of symmetries. If this can be done, then the above
theorems apply {\it mutatis mutandis}. We give here a short list of
examples.

\begin{enumerate}
\item Consider $M = S^{n+1}$ the unit $(n+1)$-dimensional sphere with
the standard metric and $N = S^n (r)$ the meridian at height
$\sqrt{1-r^2}$. The hypersurface $N$ has constant mean curvature and
is volume-degenerate since there are nontrivial Jacobi fields $w_i
(x) = x \cdot e_i$, for $i=1, \ldots ,n$ coming from the action of
the orthogonal group. Here $e_1, \ldots, e_{n+1}$ is an orthonormal
basis of ${\R}^{n+1}$. However, one may work with hypersurfaces and
functions which are invariant under the action of the $n$ hyperplanar
symmetries
\[
I_i: (x_1, \ldots, x_i, \ldots, x_{n+1}) \longrightarrow (x_1,
\ldots, -x_i, \ldots, x_{n+1})
\]
for $i=1, \ldots, n$. Namely, hypersurfaces $\tilde N\subset S^{n+1}$
such that $I_i (\tilde N) =\tilde N$ and functions $u :
S^{n+1}\longrightarrow {\mathbb R}$ such that $ u \circ I_i = u$.
Since none of the Jacobi fields is invariant under all the symmetries
$I_i$, our construction applies and the conclusion of Theorem~2 is
still valid. Moreover, when $r=1$, the equator is a minimal
hypersurface and  Theorem~1 is also valid.

\item Consider $M = B^{n+1}$ the unit ball of ${\R}^{n+1}$ endowed
with the induced metric and $N$ is a spherical cap. This example can
be dealt like the previous one and the result of Theorem~2 holds.
Moreover, when $N$ is the horizontal hyperplane, Theorem~1 also holds.

\item Consider a flat torus $T^{n+1}$ and $N$ is a pair of parallel
hyperplanes. For the sake of simplicity, assume that $T^{n+1} =
{\mathbb R}^{n+1}/ {\mathbb Z}^{n+1}$ and is identified with
$[-\frac{1}{2},\frac{1}{2}]^{n+1}$. Finally, assume that $N$ is given
by $x_{n+1} = \pm \alpha$ for some fixed $\alpha \in (0,1)$. The
action of the orthogonal group and the action of translations induce
nontrivial Jacobi fields. Again, one may work with hypersurfaces and
functions which are invariant under the action of the $n+1$
hyperplanar symmetries
\[
I_i: (x_1, \ldots, x_i, \ldots, x_{n+1}) \longrightarrow (x_1,
\ldots, -x_i, \ldots, x_{n+1})
\]
for $i=1, \ldots, n+1$ and reduce to a volume-nondegenerate problem
to show that the conclusion of Theorem~1 and Theorem~2 are valid.
\end{enumerate}

\section{Comments}\setcounter{equation}{0}

We state here a number of comments, open problems and directions for further investigations~:

\begin{enumerate}

\item In \cite{Mal}, A. Malchiodi and M. Montenegro have constructed solutions of
\begin{equation}
\e^2 \, \Delta u - u + u^p = 0,
\label{eq:neweq}
\end{equation}
which are defined on a $2$ dimensional domain and which have $0$
Neumann boundary condition. These solutions have the property that
they concentrate along the boundary of the domain and they can be
obtained for $\e$ belonging to some sequence of intervals which
converge to $0$. Behind this result, lies a very interesting
bifurcation phenomena which somehow prevents the construction to work
for all $\e$ close enough to $0$. The analysis of A. Malchiodi and M. 
Montenegro relies on the precise estimate of the least eigenvalue of 
the linearized operator about an approximate solution, and this forces them to work in $2$ dimensional domains. It would be very 
interesting to construct solutions of (\ref{eq:neweq}) which 
concentrate along a geodesic of a two dimensional manifold and also to 
extend their result to higher dimensional domains.

\item As already mentionned, M. Kowalczyk \cite{Kow} has obtained a similar result when $M$ is a two dimensional domain of ${\R}^2$ and $N$ is a line segment. Our result can be also be compared to the result of C.H. Taubes \cite{taubes} where solutions of the Seiberg-Witten equation concentrating along holomorphic curves are constructed, though our
analysis is completely different. Very closely related to our result, is the result of S. Brendle \cite{Bre-4} on the construction of solutions of the Ginzburg-Landau equation which concentrate along codimension $2$ minimal submanifolds. In all these results, solutions to nonlinear partial differential equation which concentrate along smooth submanifolds are constructed. The concentration set is always a minimal submanifold, which has to be assumed to be nondegenerate. Hence the construction holds for a generic choice of the background metric.

\item Let us also mention the construction of harmonic maps which concentrate along codimension $2$ minimal submanifolds by T. Rivi\`ere and the first author \cite{Pac-Riv} and also the work of S. Brendle concerning Yang-Mills connections in higher dimensions \cite{Bre-3}. However, in these two results, the construction does not hold for a generic choice of the metric but rather for a fairly restricted set of metrics. Finaly, the work of S. Brendle and G. Tian \cite{Bre-1} and \cite{Bre-2} should also be cited, though in these works, the difficulty is of a different nature.

\item We have not studied the case where the hypersurfaces are
singular. For example, it is known that stable minimal cones do exist
in dimension $n+1 \geq 8$ and it would be very interesting to develop
the corresponding analysis in this case.

\item  We have not studied the case where there is a nontrivial group
of isometries acting on $M$ and where the problem cannot be reduces
to a nondegenerate problem by a working in the space of functions and
hypersurfaces which are equivariant with respect to some finite group
of symmetries. For example, one may consider the case where $M$ is
the $(n+1)$-sphere with the standard metric and $N$ is the equator.
Now, let us perturb slightly the metric on $S^{n+1}$ in some
neighborhood of the north pole. Unless the perturbed metric is
invariant under the action of $I_i$, for $i=1, \ldots, n$, our method
does not apply in this setting.
\end{enumerate}

\section{The canonical profile}\setcounter{equation}{0}

In this section, we consider the case where $M={\R}$ and $\lambda
=0$. In this case, the equation (\ref{eq:2.3}) reduces to the
following second order ordinary differential equation
\begin{equation}
    \del_t^2 u - \frac{1}{2}  \, W'(u) =0 .
\label{eq:5.1}
\end{equation}
Observe that
\[
H(u,\del_t u) : = (\del_t u)^2 - W(u),
\]
is constant along solutions of (\ref{eq:5.1}). Using this property it
is easy to check that there exists a solution of (\ref{eq:5.1}),
which will be denoted by $u_\star$ in the remaining of the paper, and
which satisfies
\[
\lim_{t\rightarrow \pm \infty} u_\star (t) = \pm 1,  \qquad
\mbox{and} \qquad  u_\star (0) = 0.
\]
This solution corresponds to $H(u, \del_t u) \equiv 0$ and is
implicitly defined by
\[
t = \int_0^{u_\star(t)}  \frac{dx}{\sqrt{ W(x)}}.
\]

We define the indicial roots $\gamma_\pm >0$ by
\begin{equation}
\gamma_\pm^2 : = \frac{1}{2}\,  W''(\pm 1),
\label{eq:5.2}
\end{equation}
(observe that $\pm 1$ are minimizers of $W$ and are assumed to be
nondegenerate, hence $W''(\pm 1) >0$ and this implies that
$\gamma_\pm$ are well defined). The asymptotics of the function
$u_\star$ as $t$ tends to $\pm \infty$ are easy to derive by
linearizing (\ref{eq:5.1}) about $u\equiv \pm 1$. We find that, for
all $k \in {\mathbb N}$, there exists $c_k >0$ such that
\begin{equation}
|\del_t^k (u_\star (t) +  1) |\leq c_k \, e^{\gamma_- t} \qquad
\mbox{for all} \qquad t \leq 0,
\label{eq:5.3}
\end{equation} and
\begin{equation}
|\del_t^k ( u_\star (t) -  1) |\leq c_k \, e^{-\gamma_+ t} \qquad
\mbox{for all} \qquad t \geq 0 .
\label{eq:5.4}
\end{equation}

\section{Injectivity results}\setcounter{equation}{0}

We prove some injectivity results for ordinary differential operators
and partial differential operators whose potential is defined using
the function $u_\star$.

\subsection{Preliminary results}

For all $\zeta \in {\mathbb R}$, we set
\[
L_\zeta : =  - \, \del_t^2 + \zeta + \frac{1}{2} \, W''(u_\star) .
\]
All the injectivity results rely on the fact that, when $\zeta =0$,
the function
\[
w_\star : = \del_t u_\star ,
\]
is a bounded positive solution of the homogeneous problem
\[
L_0 \, w_\star = 0 .
\]
Furthermore $w_\star$ decays exponentially at both $+\infty$ and
$-\infty$. We introduce, for $\zeta \geq 0$, the indicial roots of
the operator $L_\zeta$ at $\pm \infty$ by
\begin{equation}
\gamma_\pm (\zeta) : =  \sqrt{\zeta + \gamma_\pm} ,
\label{eq:6.2}
\end{equation}
where $\gamma_\pm$ have been defined in (\ref{eq:5.2}). These
indicial roots are related to the asymptotic expansion near $\pm
\infty$ of the solutions of the homogeneous problem $L_\zeta \, w=0$.
For example, it follows from Cauchy's existence result for solutions
of ordinary differential equations that there exist $\underline w$
and $\overline w$ solutions of $L_\zeta w =0$, which satisfy
\[
\lim_{t\rightarrow +\infty} e^{-\gamma_+ t} \, {\overline w} (t) = 1 ,
\]
and
\[
\lim_{t\rightarrow +\infty} e^{\gamma_+ t} \, {\underline w} (t) = 1 .
\]

Our first injectivity result reads~:
\begin{lemma}
Assume that $\zeta \geq 0$ and let $w$ be a solution of $L_\zeta \, w
=0$ which is defined on $(t_1, t_2)$. Further assume that $w(t_i)=0$,
for $i=1,2$. Then $w \equiv 0$.
\label{le:fir}
\end{lemma}
{\bf Proof~:} We argue by contradiction. Given $\eta \in {\mathbb
R}$,  observe that
\begin{equation}
L_\zeta (w_\star + \eta \, w) = \zeta \, w_\star .
\label{eq:a11}
\end{equation}
Now, since we have assumed that $w$ is not identically equal to $0$,
one can choose $\eta \in {\mathbb R}$ such that the infimum of the
function $W : = w_\star + \eta \, w $, over $[t_1, t_2]$ is equal to
$0$. Observe that, since $w(t_1)=w(t_2)=0$, the infimum of $W$ is
achieved at some point $t_0 \in (t_1, t_2)$ and $L_\zeta W \leq 0$ at
this point. When $\zeta >0$, this last inequality clearly contradicts
(\ref{eq:a11}). When $\zeta =0$, observe that $W$ is a solution of
some second order linear ordinary differential equation and that
$W(t_0) = \del_t W(t_0) = 0$, hence $W\equiv 0$. This contradicts the
fact that $W(t_i)\neq 0$ for $i=1,2$. \hfill $\Box$

Our second injectivity result classifies the set of $\zeta$ for which
there exists a bounded solution of $L_\zeta w = 0$, which is defined
on ${\mathbb R}$ or on a half line. Given $\delta_\pm \in {\R}$, we
define $\delta := (\delta_-, \delta_+)$ and the function
\begin{equation}
\varphi_\delta (t) : =  (1+ e^t)^{\delta_+} \, (1+e^{-t})^{\delta_-}.
\label{eq:weight}
\end{equation}
In particular, $\varphi_\delta(t) \sim e^{\delta_+ t}$ at $+\infty$
and $\varphi_\delta(t) \sim e^{-\delta_- t}$ at $-\infty$. This
definition being understood, we now prove the~:
\begin{lemma}
Assume that $\zeta \neq 0$. Let $w$ be a solution of $L_\zeta \, w =
0$ which is defined on ${\mathbb R}$ (or on $(-\infty, t_0)$ or
$(t_0, +\infty)$, in which case we ask that $w(t_0)=0$). Further
assume that $|w|$ is bounded by a constant times the function
$\varphi_\delta$ for  $\delta_\pm < \gamma_\pm (\zeta)$. Then $w
\equiv 0$.
\label{le:sir}
\end{lemma}
{\bf Proof~:}
The proof of this result is almost identical to the proof of the
previous result. The key observation being that, under the above
assumptions, any solution of $L_\zeta w =0$ defined on a half line
decays faster than $w_\star$ at infinity. This follows at once from
the fact that the indicial roots $\gamma_\pm (\zeta) $ of $L_\zeta$
are larger than the indicial roots $\gamma_\pm$ of $L_0$.

For example, assume that $w$ is defined on $(t_0, +\infty)$. Any
solution of the homogeneous problem $L_\zeta w=0$ is a linear
combination of $\overline w$, a solution which blows up exponentially
at $+\infty$ like  $e^{\gamma_+ (\zeta) t}$, and $\underline w$, the
solution which decays exponentially at $+\infty$ like
$e^{-\gamma_+(\zeta) t}$. Since $|w|$ is bounded by $e^{\delta_+t}$
for some $\delta_+ < \gamma_+ (\zeta)$, we conclude that $w$ is
collinear to $\underline w$. Now, $\gamma_+ (\zeta) > \gamma_+$ hence
$w$ decays faster than $w_\star$ at $+\infty$. Once this is known,
the proof of the result reduces to the proof of Lemma~\ref{le:fir}.
\hfill $\Box$

The set of solutions of $L_\zeta w = 0$ is two dimensional and there
exists a unique $w^-_\zeta$ solution of $L_\zeta w_\zeta^- =0$ which
is defined on all $\R$ and which satisfies
\[
\lim_{t\rightarrow -\infty} e^{- \gamma_- (\zeta)  \, t}\, w_\zeta^-  (t) = 1.
\]
As already mentioned, this essentially follows from Cauchy's
existence result for solutions of ordinary differential equations.
When $\zeta \neq 0$, the previous Lemma implies that the function
$w_\zeta^-$ does not vanish and furthermore blows up exponentially at
$+\infty$ like $t \rightarrow e^{\gamma_+ (\zeta) \, t}$. Similarly,
there exists a unique $w^+_\zeta$ solution of $L_\zeta w_\zeta^+ =0$
on $\R$ which satisfies
\[
\lim_{t\rightarrow +\infty} e^{\gamma_+ (\zeta)  \, t} \, w_\zeta^+  (t) = 1.
\]
Again, when $\zeta \neq 0$, this function does not vanish and blows
up exponentially at $-\infty$ like $t \rightarrow e^{- \gamma_-
(\zeta) \, t}$.

\subsection{Injectivity results}

Assume that $(N,h)$ is a compact $n$-dimensional Riemannian manifold
with or without boundary. Further assume that $(N,h)$ is at least
${\cal C}^{1, \alpha}$. This means that one can choose local
coordinate charts on $N$ in which the coefficients of the metric $h$
are ${\cal C}^{1, \alpha}$ functions. We define on the product space
${\mathbb R}\times N$ the partial differential operator
\begin{equation}
{\mathfrak L}_h : =  - \del_t^2 - \Delta_h + \frac{1}{2} \, W''(u_\star) ,
\label{eq:8.1}
\end{equation}
where $\Delta_h$ is the Laplace-Beltrami operator on $(N,h)$. Using
the result of Lemma \ref{le:fir}, we get the~:
\begin{corollary}
Assume that $w$ is a solution of ${\mathfrak L}_h \, w =0$ which is
defined on $(t_1, t_2) \times N$. Further assume that $w=0$ on
$\{t_i\}\times N$, for $i=1,2$ and that $w$ has $0$ Neumann boundary
data on $(t_1, t_2) \times \del N$ if $\del N$ is not empty. Then $w
\equiv 0$.
\label{co:7.1}
\end{corollary}
{\bf Proof~:} We denote by $(\phi_j, \lambda_j)_{j \geq 0}$ the
eigendata of $\Delta_h$ (with Neumann boundary conditions when the
boundary of $N$ is not empty). Namely
\[
\Delta_h \phi_j =  - \lambda_j \, \phi_j ,
\]
with $\lambda_i \leq \lambda_{i+1}$. We also assume that the
eigenfunctions are normalized so the their $L^2$ norm on $N$ is $1$.
We decompose the function $w$ defined on $(t_1,t_2) \times N$ as
\[
w(t,y) = \sum_{j \in {\mathbb N}} w_j(t) \, \phi_j (y) .
\]
Then $w_j$ is a solution of $L_{\lambda_j} w_j =0$ and the result of
Lemma~\ref{le:fir} implies that $w_j\equiv 0$. This completes the
proof.  \hfill $\Box$

Using similar arguments together with Lemma~\ref{le:sir}, we also get the~:
\begin{corollary}
Assume that $w$ is a solution of ${\mathfrak L}_h \, w =0$ which is
defined on ${\R}\times N$. Further assume that $w$ is bounded by a
constant times $\varphi_\delta$ for some $\delta_\pm < \gamma_\pm$.
Then $w$ only depends on $t$ and there exists a constant $c
\in{\mathbb R}$ such that $w = c \, w_\star$.
\label{co:frc}
\end{corollary}
Before we proceed, let us observe that the above result holds under a
slightly more general assumption, namely that $\delta_\pm <
\gamma_\pm (\lambda_1)$, where $\lambda_1$ is the first nonzero
eigenvalue of the Laplace-Beltrami operator $\Delta_h$.

{\bf Proof~:} Again, we decompose the function $w$ defined on
${\mathbb R} \times N$ as
\[
w(t,y) = \sum_{j \in {\mathbb N}} w_j(t) \, \phi_j (y)
\]
Then $w_j$ is a solution of $L_{\lambda_j} w_j =0$ and the result of
Lemma~\ref{le:sir} implies that $w_j\equiv 0$, for all $j \neq 0$.
When $j=0$, all bounded solutions of $L_0 w =0$ have to be collinear
to $w_\star$. This completes the proof. \hfill $\Box$

Let $\Delta$ denote the Laplacian in ${\mathbb R}^n$ endowed with the
Euclidean metric. We define the elliptic operator
\[
{\mathfrak L} : =  - \del_t^2 - \Delta  + \frac{1}{2}\, W''(u_\star) .
\]
The result of Lemma~\ref{le:sir} also implies the~:
\begin{corollary}
Assume that $w \in L^{\infty} ({\R}\times {\R}^n)$ is a solution of
${\mathfrak L} \, w =0$. Then $w$ only depends on $t$ and there
exists a constant $c \in{\mathbb R}$ such that $w = c \, w_\star$.
\label{le:6.3}
\end{corollary}
This result seems to be standard and appear, without a proof in
\cite{Maz-Pac}, but we have not been able to find a precise reference
for it. Since this is a key result for the forthcoming argument, we
give here a detailed proof.

{\bf Proof~:} We denote by ${\cal S}({\mathbb R}^k)$, the space of
smooth rapidly decaying functions which are defined on ${\R}^k$. This
space is endowed with the family of semi norms
\[
[\, \phi \, ]_{k,l} : =  || (1+|z|^k) \, \nabla^l \phi||_{L^\infty} ,
\]
for all $k, l \in {\mathbb N}$, where $z$ denotes the variable in
${\mathbb R}^k$. The dual space ${\cal S}'({\mathbb R}^k)$ is the
space of tempered distributions \cite{Sch}.

Let ${\cal F}$ denote the Fourier transform in ${\R}^n$ and
$\overline{\cal F}$ the Fourier inverse transform. We define, for all
$\phi \in {\cal S}({\R}\times {\R}^n)$,
\[
T ( \phi ) : =  \int_{\R} \, \langle w(t, \cdot), \overline{\cal F}
(\phi (t, \cdot)) \rangle_{{\cal S}',{\cal S}} \, dt.
\]
This clearly defines a tempered distribution $T \in {\cal
S}'({\R}\times {\R}^n)$. Let us denote by $\xi \in {\R}^n$ the dual
variable of $z \in {\R}^n$. Using the fact that $w$ is a solution of
${\mathfrak L} \, w =0$, we get
\begin{equation}
\langle L_{|\xi|^2}  \, T, \Psi \rangle_{{\cal S}', {\cal S}}: =
\langle T, L_{|\xi|^2} \Psi \rangle_{{\cal S}', {\cal S}} = 0,
\label{eq:6.4}
\end{equation}
for any smooth function $\Psi \in {\cal S} ({\R}\times {\R}^n)$.

We claim that the support of $T$ is included in ${\R}\times \{0\}$.
Indeed, choose any smooth function  $\psi : {\R}\times {\R}^n
\longrightarrow {\R}$ with compact support in ${\R}\times ({\R}^n
-\{0\})$. We define
\[
\Psi (t, \xi) : =  \frac{1}{\alpha (\xi)} \left( w^+_{|\xi|^2}
(t)\int_{t}^{+\infty} w^-_{|\xi|^2} (s)\psi (s, \xi) ds +
w^-_{|\xi|^2} (t)  \int_{-\infty}^t w^+_{|\xi|^2} (s) \psi (s, \xi)ds
\right) ,
\]
where
\[
\alpha (\xi) : = w^-_{|\xi|^2} (t)\, \del_t w^+_{|\xi|^2} (t) -
\del_t w^-_{|\xi|^2} (t) \, w^+_{|\xi|^2} (t),
\]
is the Wronskian of the two independent solutions $w^\pm_{|\xi|^2}$
of the homogeneous problem $L_{|\xi|^2} w=0$ which have been defined
at the end of \S 7.1 (hence $\alpha$ does not vanish and does not
depend on $t$ !). We claim that $L_{|\xi|^2} \Psi = \psi$ and also
that $\Psi \in {\cal S}({\R}\times {\R}^n)$. The first claim follows
at once from the fact that $w^\pm_{|\xi|^2}$ are solutions of
$L_{|\xi|^2} w=0$. For the second claim, observe that the function
$\psi$ has compact support in ${\R}\times ({\R}^n - \{0\})$, hence
$\Psi (t,\xi) \equiv 0$ for all $|\xi|$ large enough (say $|\xi |
\geq c$) and $|\xi|$ small enough (say $|\xi|\leq 1/c $). To show
that $\Psi$ is rapidly decaying in $t$ when $1/c \leq |\xi| \leq c$,
we use the fact that, for $\xi \neq 0$, the function $w^+_{|\xi|^2}$
is exponentially decaying at $+\infty$ and the function
$w^-_{|\xi|^2}$ is exponentially decaying at $-\infty$. This implies
at once that $\Psi \in {\cal S} ({\R}\times {\R}^n)$. Therefore, we
conclude that
\[
\langle T, \psi \rangle_{{\cal S}', {\cal S}} = 0,
\]
for all $\psi : {\R}\times {\R}^n \longrightarrow {\R}$ with compact
support in ${\R}\times ({\R}^n -\{0\})$. This proves the claim.

By a classical result in the theory of distributions \cite{Sch}, we
know that $T$ is the linear combination of derivatives (of bounded
order), with respect to $t$ and $\xi_j$, $j=1, \ldots, n$, of
measures with support on ${\R}\times \{0\}$. Performing the Fourier
transform backward in the $\xi$ variable, we see that the function
$w(t, \cdot)$ depends polynomially on the coordinates $z_j$ of $z\in
{\R}^n$. This, together with the fact that $w$ is bounded in the $z$
variable, implies that the function $w$ only depends on $t$ and hence
$w = c \, w_\star$. This completes the proof of the result. \hfill
$\Box$

Given $\gamma \in {\R}$, we define the elliptic operator
\[
{\mathfrak L}_\gamma : =  - \del_t^2 - \Delta  + \gamma .
\]
Following the proof of Lemma~\ref{le:6.3}, we have~:
\begin{lemma}
Assume that $\gamma  > 0$ and that $u \in L^{\infty} ({\R}\times
{\R}^n)$ is a solution of ${\mathfrak L}_\gamma \, w =0$. Then
$w\equiv 0$.
\label{le:6.4}
\end{lemma}

\section{Mapping properties of a model operator} \setcounter{equation}{0}

In this section we study the mapping properties of the operator
${\mathfrak L}_h$, given in (\ref{eq:8.1}), when this operator is
defined between weighted H\"older spaces.

\subsection{Function spaces}

Assume that $(N,h)$ is a compact $n$-dimensional compact Riemannian
manifold with or without boundary. To begin with, let us define the
weighted spaces we will work with~:
\begin{definition}
Given $\ell \in {\N}$, $\alpha \in (0,1)$ and $\delta := (\delta_-,
\delta_+) \in {\R}^2$, we define the weighted space ${\cal C}^{\ell,
\alpha}_\delta ({\R}\times N)$ to be the space of functions which are
$\ell$ times differentiable, whose $\ell$-th partial derivatives are
H\"older of exponent $\alpha$ and for which  the weighted norm
\[
\|u\|_{{\cal C}^{\ell, \alpha}_\delta ({\R}\times  N)} : = \|
\varphi_{-\delta} \, u\|_{{\cal C}^{\ell, \alpha} ({\R}\times N)},
\]
is finite. Here by definition
\[
\|u\|_{{\cal C}^{\ell, \alpha} ({\R}\times N)} : =  \sum_{j=0}^\ell
\|  \nabla^j  u \|_{L^\infty({\R}\times N)}+ \sup_{p \neq q \, \in \,
{\R}\times N}  \, \frac{|\nabla^\ell u  (p)- \nabla^\ell u
(q)|}{d(p,q)^\alpha},
\]
is the standard (unweighted) H\"older norm on ${\R}\times N$ and $d$
is the geodesic distance in ${\R}\times N$ for the product metric.
\end{definition}
Roughly speaking, this is the space of functions which, together with
their partial derivatives, are bounded by $e^{\delta_+ t}$ on $(0,
+\infty)\times N$ and are bounded by $e^{-\delta_- t}$ on $(-\infty,
0)\times N$.

We finally define a $1$-codimensional subspace of ${\cal C}^{\ell,
\alpha}_\delta ({\R}\times  N)$.
\begin{definition}
Given $\ell \in {\N}$,  $\alpha \in (0,1)$ and $\delta := (\delta_-,
\delta_+) \in {\R}^2$ such that $\delta_\pm < \gamma_{\pm}$. The
space ${\cal D}^{\ell, \alpha}_\delta ({\R}\times N)$ is defined to
be the closed subspace of functions $u \in {\cal C}^{\ell,
\alpha}_\delta ({\R}\times N)$ which are $L^2$ orthogonal to
$w_\star$, i.e.
\begin{equation}
\int_{{\R}\times N} u(t,y) \, w_\star (t)  \, dt \, d a_h =0 .
\label{eq:8.2}
\end{equation}
Naturally, this space is endowed with the induced norm.
\end{definition}
The restriction $\delta_\pm < \gamma_\pm$ is needed to ensure the
convergence of this integral in (\ref{eq:8.2}), i.e. that $u \,
w_\star \in L^1 ({\mathbb R}\times N)$.

In the case where $N$ has a boundary, we define, for $k \geq 1$
\[
\left[ {\cal D}^{\ell, \alpha}_\delta ({\R}\times N)\right]_{0} ,
\]
to be the subspace of functions of  ${\cal D}^{\ell, \alpha}_\delta
({\R}\times N)$ which have $0$ Neumann boundary condition on
${\R}\times \del N$. In the subsequent sections, it will be
convenient to adopt the notation $\left[ {\cal D}^{\ell,
\alpha}_\delta ({\R}\times N)\right]_{0}$ for ${\cal D}^{\ell,
\alpha}_\delta ({\R}\times N)$ even when $N$ has no boundary, being
understood that the condition on the boundary data is void in this
latter case.

\subsection{Mapping properties}

Recall that we have defined on the product space ${\mathbb R}\times
N$ the partial differential operator
\begin{equation}
{\mathfrak L}_h : =  - \del_t^2 - \Delta_h  + \frac{1}{2} \, W''(u_\star) ,
\label{eq:8.1b}
\end{equation}
where $\Delta_h$ is the Laplace-Beltrami operator on $(N,h)$. We now
assume that $(N, h)$ is at least ${\cal C}^{\ell-1, \alpha}(N)$ for
some $\ell \geq 2$. Recall that this means that there exists local
coordinate $y_1, \ldots, y_n$ on $N$ for which the coefficients
$h_{ij}$ of the metric
\[
h : =  \sum_{i,j} h_{ij}\, dy_i \otimes dy_j ,
\]
are ${\cal C}^{\ell-1, \alpha}$ functions on $N$. Clearly, the operator
\begin{equation}
{\mathfrak L}_h : \left[ {\cal D}^{\ell, \alpha}_\delta ({\R}\times
N)\right]_{0} \longrightarrow  {\cal D}^{\ell-2, \alpha}_\delta
({\R}\times N) ,
\label{eq:8.20}
\end{equation}
is well defined and bounded, for any $\delta \in {\R}^2$. It is well
known that the mapping properties of the above defined operator
depends crucially on the choice of the weight parameter $\delta$.
Indeed, the set of indicial roots of ${\mathfrak L}_h$ at $+ \infty$
(resp. $-\infty$)  is defined by ${\cal I}_+$ (resp. ${\cal I}_-$)
where
\[
{\cal I}_\pm : = \{\gamma_\pm (\lambda_j) \, : \, j \geq 0\},
\]
and $\lambda_j$ are the eigenvalues of $\Delta_h$ on $N$. Now, if
$\delta_\pm \notin {\cal I}_\pm$, the operator (\ref{eq:8.20}) can be
proven to have closed range and to be Fredholm. While, if $\delta_-$
or $\delta_+$ is an indicial root, then the operator (\ref{eq:8.20})
does not have close range, and hence is not Fredholm.

The result of Corollary~\ref{co:frc} yields the injectivity of the
operator ${\mathfrak L}_h$ when $\delta_\pm < \gamma_\pm$. This,
together with a "duality argument" implies that the operator
(\ref{eq:8.20}) is surjective provided $\delta_\pm > - \gamma_\pm$
and $\delta_\pm \notin {\cal I}_\pm$ (the duality argument does hold
{\it stricto sensus} when the operator is defined between weighted
Lebesgue spaces and the corresponding result in weighted H\"older
spaces is then obtained through elliptic regularity theory). In
particular, the operator (\ref{eq:8.20}) is an isomorphism if the
weight $\delta_\pm \in (-\gamma_\pm, \gamma_\pm)$. Most of these
mapping properties of ${\mathfrak L}_h$ will not be needed. Indeed,
we will only need the later claim on the range of weights for which
the operator is an isomorphism. Therefore, in the next Proposition,
we  concentrate on the proof of this fact.

In addition, we will also show that the inverse of ${\mathfrak L}_h$
is  bounded independently of the choice of the metric $h$ on $N$,
provided an uniform ellipticity condition is fulfilled. Hence, we now
assume that $\ell \geq 2$ and that there exists a constant
$\Lambda >0$ such that, for all $\xi \in {\R}^n$ and all $y \in N$
\begin{equation}
\sum_{i,j} h^{ij} (y) \xi_i\, \xi_j \geq \Lambda \, |\xi|^2 ,
\label{eq:8.22}
\end{equation}
and
\begin{equation}
\| h_{ij} \|_{{\cal C}^{\ell-1, \alpha} (N)} \leq \frac{1}{\Lambda}.
\label{eq:8.222}
\end{equation}
The first condition will ensure that the Laplace-Beltrami operator on
$N$ is uniformly elliptic while the second condition will ensure
uniform H\"older estimates, independent of the choice of the metric.

We prove the following~:
\begin{proposition}
Assume that $\ell \geq 2$ is fixed and that $h$ satisfies
(\ref{eq:8.22}) and (\ref{eq:8.222}). Further assume that $
\delta_{\pm} \in (- \gamma_\pm ,  \gamma_\pm)$. Then, the operator
${\mathfrak L}_h$ defined in (\ref{eq:8.20}) is an isomorphism and
there exists a constant $c >0$, only depending on $\Lambda$, such
that, for all $w \in [{\cal D}^{\ell, \alpha}_\delta ({\R}\times
N)]_{0}$ we have
\[
||w||_{{\cal C}^{\ell, \alpha}_\delta ({\R}\times N)}\leq c \,
||{\mathfrak L}_h \, w||_{{\cal C}^{\ell-2, \alpha}_\delta
({\R}\times N)}.
\]
\label{pr:8.1}
\end{proposition}
{\bf Proof :} As already mentioned, the injectivity of the operator
${\mathfrak L}_h$ follows at once from Corollary~\ref{co:frc}.
Indeed, given the range in which we have chosen $\delta_\pm$, any
solution of ${\mathfrak L}_h w=0$ which is bounded by
$\varphi_\delta$ has to be collinear to $w_\star$. Since functions in
$[{\cal D}^{\ell, \alpha}_\delta({\R}\times N)]_{0}$ are orthogonal
to $w_\star$ in $L^2 ({\mathbb R}\times N)$, we conclude that $w
\equiv 0$.

We now prove that the operator ${\mathfrak L}_h$ defined in
(\ref{eq:8.20}) is surjective.  To this aim, we decompose any
function $f \in {\cal D}^{\ell-2, \alpha}_\delta ({\R}\times N)$ as
\[
f(t,y) = f_0(t) + \tilde f(t,y) ,
\]
where, for each $t$, the function $\tilde f(t, \cdot)$ is orthogonal
to the constant function $1$ in $L^2 (N)$. The proof is now
decomposed into 3 steps.

\noindent
{\bf Step 1.} We define
\[
w_0 (t) : = w_\star (t)  \, \left( \alpha_0  + \int_0^t \,
w_\star^{-2}(s) \, \int_{-\infty}^s w_\star(r) \, f_0(r) \, dr \, ds
\right) ,
\]
where the constant $\alpha_0 \in {\R}$ is chosen so that
\[
\int_{\R} w_0 \, w_\star \, dt = 0 .
\]
Using the fact that $w_\star \sim e^{\gamma_- t}$ at $-\infty$,
together with $-\gamma_- < \delta_-< \gamma_-$, we conclude that
\[
\sup_{(-\infty,0]} e^{\delta_- t} \, |w_0|\leq c_0 \, \sup_{\R}
e^{\delta_- t}\, |f_0|.
\]
Since$ f \in {\cal D}^{\ell-2, \alpha}_\delta ({\R}\times N)$, we have
\[
\int_{\R} f_0 \, w_\star \, dt = 0 .
\]
Therefore, we have the alternative definition of $w_0$ as
\[
w_0 (t) : = w_\star (t)  \, \left( \alpha_0  - \int_0^t \,
w_\star^{-2}(s) \, \int_s^{+\infty} w_\star(r) \, f_0(r) \, dr \, ds
\right) ,
\]
and, using this time the fact that $w_\star \sim e^{- \gamma_+ t}$ at
$+\infty$ together with $-\gamma_+ < \delta_+ < \gamma_+$, we
conclude that
\[
\sup_{\R} e^{-\delta_+ t} \, |w_0|\leq c_0 \, \sup_{\R} e^{-\delta_+
t} \, |f_0|,
\]
for some constant $c_0 >0$ which does not depend on $f_0$.

\noindent
{\bf Step 2.} Making use of the result of Corollary~\ref{co:7.1}, we
can solve, for each $T >0$
\[
{\mathfrak L}_h \tilde w_T = \tilde f, \qquad \mbox{in} \qquad (-T,T) \times N,
\]
with $\tilde w_T = 0$ on $\{\pm T\} \times N$ and $0$ Neumann
boundary data on $(-T,T) \times \del N$ if $N$ has a boundary. The
superscript $\sim$ is meant to recall that we are working with
functions $(t,y) \rightarrow \tilde g(t,y)$ for which the function
$\tilde g(t, \cdot)$ is for each $t$ orthogonal to the constant
function $1$ in $L^2 (N)$.

\noindent
{\bf Step 3.} We claim that
\[
|| \varphi_{-\delta} \, \tilde w_T||_{L^\infty ((-T,T)\times N)}\leq
c \, ||  \varphi_{-\delta} \, \tilde f||_{L^\infty ({\R}\times N)},
\]
for some constant $c >0$ which depends on $\Lambda$ but does not
depend on $T$ nor on $\tilde f$.

Indeed, choose $t_0 >0 $ large enough so that
\[
{\mathfrak p} : = \min \left( \inf_{(-\infty , -t_0)} \left(
\frac{1}{2} \, W''(u_\star)  - \delta_-^2 \right) ,  \inf_{(t_0,
+\infty)} \left(\frac{1}{2} \, W''(u_\star) \right) -\delta_+^2
\right) > 0 ,
\]
then, the potential in the operator ${\mathfrak L}_h$ is bounded from
below by ${\mathfrak p}$ in $({\R}-(-t_0, t_0)) \times N$ and hence
${\mathfrak L}_h$ satisfies the maximum principle in this set.
Moreover, we have
\[
{\mathfrak L}_h \, e^{-\delta_- t} =  \left( \frac{1}{2} \,
W''(u_\star)  - \delta_-^2 \right)  \, e^{-\delta_- t} \geq
{\mathfrak p} \, e^{-\delta_- t} ,
\]
in $(-\infty, -t_0) \times N$ and
\[
{\mathfrak L}_h \, e^{\delta_+ t} =  \left( \frac{1}{2} \,
W''(u_\star)  - \delta_+^2 \right)  \, e^{\delta_+ t} \geq {\mathfrak
p} \, e^{\delta_+ t},
\]
in $(t_0, +\infty) \times N$. Hence, the function $t \rightarrow
e^{\delta_+ t}$ (resp. $t \rightarrow e^{-\delta_- t}$) can be used
as barrier in $(t_0, +\infty)\times N$ (resp. in $(-\infty,
t_0)\times N$) to prove, for all $T > t_0$, the inequality
\begin{equation}
|| \varphi_{-\delta} \tilde w ||_{L^\infty ((-T,T) \times N)} \leq c
\left(  || \varphi_{-\delta} \tilde w ||_{L^\infty ([- t_0, t_0]
\times N)} + || \varphi_{-\delta} \, \tilde f ||_{L^\infty
((-T,T)\times N)} \right),
\label{eq:hhh}
\end{equation}
for some constant $c >0$ which only depends on $\delta$ and on $\mathfrak p$.

We can now complete the proof of the claim. We argue by contradiction
and assume that there exists a sequence of metrics $h_k$ satisfying
(\ref{eq:8.22}) and (\ref{eq:8.222}), a sequence $T_k >0$, a sequence
of functions $\tilde f_k$, and a sequence of solutions of ${\mathfrak
L}_h \tilde w_k = \tilde f_k$ in $(-T_k,T_k) \times N$, with $\tilde
w_k =0$ on $\{\pm T_k\} \times N$ and $0$ Neumann boundary data on
$(-T_k,T_k) \times \del N$ if $N$ has a boundary, such that
\[
|| \varphi_{-\delta} \, \tilde w_k||_{L^\infty ((-T_k,T_k)\times N)} = 1,
\]
and
\[
\lim_{k \rightarrow \infty} || \varphi_{-\delta} \, \tilde
f_k||_{L^\infty ({\R}\times N)} =0.
\]
Furthermore, $\tilde f(t, \cdot)$ and $\tilde w(t, \cdot)$ are, for
all $t$, orthogonal to $1$ in $L^2(N)$.

Observe that the claim is certainly true when $T_k$ remains bounded,
hence we may well assume that the sequence $T_k$ tends to $+\infty$.
Since both $\tilde w_k$ and $\tilde f_k $ are uniformly bounded on
compact subsets, we deduce from elliptic estimates that the sequence
of functions $\tilde w_k$ is uniformly bounded in ${\cal C}^{1,
\alpha}$ topology on any compact of ${\R}\times N$. Now, Ascoli's
Theorem together with a standard diagonal argument implies that, up
to a subsequence and for some $\alpha' < \alpha$, the sequence of
metrics $h_k$ converges (in ${\cal C}^{1, \alpha'} (N)$ topology) to
$h_\infty$, the sequence of functions $\tilde w_k$ converges (in
${\cal C}^{1, \alpha'} (N)$ topology) to $\tilde w_\infty$. Passing
to the limit in the equation satisfied by $\tilde w_k$, we conclude
that $\tilde w_\infty$ is a weak solution of ${\mathfrak
L}_{h_\infty} \, \tilde w_\infty = 0$ in ${\R}\times N$, which is
bounded by $\varphi_\delta$. But Corollary~\ref{co:frc}  and the
choice of the parameter $\delta_\pm$ in $(-\gamma_\pm, \gamma_\pm)$
imply that $\tilde w_\infty = c \, w_\star$, for some constant $c \in
{\R}$. Since the function $\tilde w_\infty (t, \cdot)$ is by
definition orthogonal to $1$ in $L^2(N)$, we conclude that $c=0$,
hence $\tilde w_\infty \equiv 0$.

Finally, (\ref{eq:hhh}) implies that, for each $k \in {\mathbb N}$
sufficiently large (so that $T_k \geq t_0$) we have
\[
|| \varphi_{-\delta} \tilde w_k ||_{L^\infty ([-T_k, T_k] \times N)}
= 1 \leq  c   \left(  || \varphi_{-\delta}\tilde w_k ||_{L^\infty
([-t_0, t_0] \times N)} +|| \varphi_{-\delta} \tilde f_k ||_{L^\infty
({\R}\times N)} \right).
\]
Passing to the limit as $k$ tends to $\infty$, we conclude that
\[
1 \leq  c \, ||\varphi_{-\delta} \, \tilde w_\infty ||_{L^\infty
([-t_0, t_0] \times N)},
\]
which implies that $\tilde w_\infty$ is not identically equal to $0$.
A contradiction. This ends the proof of the claim.

\noindent
{\bf Step 4.} Using the result of Step 2 and 3, standard elliptic
estimate and Ascoli's Theorem, we can pass to the limit as $T$ tends
to $+\infty$ in the sequence $\tilde w_T$ and obtain a function
$\tilde w$ solution of  ${\mathfrak L}_h \, \tilde w  = \tilde f $ in
${\R}\times N$. Furthermore, the result of Step 2 implies that
\[
|| \varphi_{-\delta}  \tilde w ||_{L^\infty ({\R}\times N)} \leq c
\, || \varphi_{-\delta}  \tilde f ||_{L^\infty ({\R}\times N)} ,
\]
for some constant $c >0$ which only depends on $\Lambda$. Collecting
this result and the result of Step 1, we conclude that $w = w_0+
\tilde w$ is a solution of ${\mathfrak L}_h \, w= f$ which satisfies
(\ref{eq:8.2}). Furthermore, using Schauder's estimates, we also
conclude that
\[
|| w||_{{\cal C}^{\ell, \alpha}_\delta ({\R}\times N)} \leq c  \, ||
f||_{{\cal C}^{\ell-2, \alpha}_\delta({\R}\times N)},
\]
for some constant $c>0$ which only depends on $\Lambda$. In
particular $w \in [{\cal D}^{\ell, \alpha}_\delta ({\R}\times
N)]_{0}$. This completes the proof of the result. \hfill $\Box$

Once Proposition~\ref{pr:8.1} is proven, it is easy to see that~:
\begin{proposition}
Assume that $\delta_\pm < \gamma_\pm$.  If $w \in [{\cal D}^{\ell,
\alpha}_\delta ({\R}\times N)]_{0}$ and $f \in {\cal
D}^{\ell-2,\alpha}_\delta ({\R}\times N)$ satisfy ${\mathfrak L}_h \,
w =f$ and if in addition
\[
\int_{\R}  \, f (t,y) \, w_\star (t) \, dt = 0,
\]
for all $y \in N$, then
\[
\int_{\R}  \, w (t,y) \, w_\star (t) \, dt = 0 ,
\]
for all $y \in N$.
\label{pr:8.2}
\end{proposition}
{\bf Proof~:} Recall that $(\phi_j, \lambda_j)_{j\geq 0}$ are the
eigendata of $\Delta_h$ (with $0$ Neumann boundary condition if $\del
N$ is not empty). For all $j \geq 1$, we multiply the equation
${\mathfrak L}_h \, w = f$ by $ w_\star \, \phi_j$ and integrate by
parts over ${\R}\times N$. We obtain
\[
\int_{{\R}\times N}  w \, {\mathfrak L}_h (\phi_j \, w_\star) \, dt
\, da_h = 0 .
\]
Since ${\mathfrak L}_h (w_\star \, \phi_j ) = \lambda_j  \, w_\star
\, \phi_j$, we conclude that
\[
\lambda_j \, \int_{N}  \phi_j \, \left( \int_{\R} w \, w_\star  \, dt
\right) \, da_h = 0 .
\]
When $j \neq 0$, $\lambda_j \neq 0$, hence this implies that the function
\[
\Phi(y) : = \int_{\R} w (t,y) \, w_\star(t)  \, dt ,
\]
is orthogonal to $\phi_j$ in $L^2 (N)$. By construction this function
is also orthogonal to $\phi_0$, which is the constant function, in
$L^2(N)$. Since the $(\phi_j)_{j\geq 0}$ form a Hilbert basis of
$L^{2}(N)$, we conclude that $\Phi \equiv 0$. \hfill $\Box$

Before we proceed further, let us comment on our choice of function
spaces. Observe that, in (\ref{eq:8.20}), we could replace the spaces
${\cal D}^{\ell, \alpha}_\delta ({\R}\times N)$ by the spaces ${\cal
C}^{\ell, \alpha}_\delta ({\R}\times N)$ and define
\[
\tilde{\mathfrak L}_h : \left[ {\cal C}^{\ell, \alpha}_\delta
({\R}\times N)\right]_{0} \longrightarrow  {\cal C}^{\ell-2,
\alpha}_\delta ({\R}\times N) .
\]
Then all the above results about the set of weights for which the
operator is Fredholm remain true. However, this time, the injectivity
of $\tilde{\mathfrak L}_h$ only holds provided $\delta_\pm < -
\gamma_\pm$ and, using a "duality argument", this implies that the
operator $\tilde{\mathfrak L}_h$ is now surjective when $\delta_\pm >
\gamma_\pm$ are not indicial roots. Hence, this choice of function
spaces would force us to work in a space of functions which blow up
exponentially at $\pm \infty$ and this would not be suitable for the
forthcoming nonlinear analysis.

\section{Mapping properties of a singularly perturbed linear
operator}\setcounter{equation}{0}

\subsection{Function spaces}

Assume that $(N, h)$ is a $n$-dimensional compact Riemannian
manifold, with or without boundary and further assume that the metric
$h$ on $N$ satisfies (\ref{eq:8.22}) and (\ref{eq:8.22}), for some
fixed $\ell \geq 2$ and some fixed constant $\Lambda >0$. We now turn
to the study of the operator
\begin{equation}
{\mathfrak L}_\e : = - \e^2 \, \left( \del_t^2 + \Delta_h \right) +
\frac{1}{2} \, W''(u_\star (\cdot/\e)),
\label{eq:9.1}
\end{equation}
where $\e \in (0,1)$ is a fixed parameter. This operator still
depends on $h$ but, since we now focus our attention on its
dependence with respect to the parameter $\e$, we omit the subscript
$h$. Taking the parameter $\e$ into account, we now define~:
\begin{definition}
Given, $\ell \in {\N}$, $\alpha \in (0,1)$, $\delta : =  (\delta_-,
\delta_+) \in {\R}^2$ and $\e \in (0,1)$, we define the weighted
space ${\cal C}^{\ell , \alpha}_{\delta,\e} ({\R}\times N)$ to be the
space of functions which are $\ell$ times differentiable, whose
$\ell$-th partial differential are H\"older of exponent $\alpha$ and
for which the weighted norm
\[
|| u ||_{{\cal C}^{\ell, \alpha}_{\delta,\e} ({\R}\times N)} : =
||\varphi_{-\delta} (\cdot/\e) \, u ||_{{\cal C}^{\ell, \alpha}_{\e}
({\R}\times N)},
\]
is finite. Here, by definition
\[
|| u ||_{{\cal C}^{\ell, \alpha}_{\e} ({\R}\times N)} : =
\sum_{j=0}^\ell  \e^j \, \|\nabla^j u \|_{L^\infty ({\R}\times N)}+
\e^{\ell + \alpha} \, ||u||_{{\cal C}^{0, \alpha}({\R}\times N)} .
\]
\end{definition}
When $\delta_\pm < \gamma_\pm$, we also define the spaces ${\cal
D}^{\ell, \alpha}_{\delta ,\e} ({\R}\times N)$ as the
$1$-codimensional subspace of functions $u \in {\cal C}^{\ell,
\alpha}_{\delta ,\e} ({\R}\times N)$ which are $L^2$ orthogonal  to
the function $w_\star (\cdot /\e)$. For fixed $\e$ the result \S 8
apply {\it verbatim} and this shows that
\[
{\mathfrak L}_\e : [{\cal D}^{\ell, \alpha}_{\delta ,\e} ({\R}\times
N)]_{0} \longrightarrow {\cal D}^{\ell -2, \alpha}_{\delta, \e}
({\R}\times N) ,
\]
is an isomorphism, provided $\delta_\pm \in (-\gamma_\pm,
\gamma_\pm)$. Our aim is now to understand the mapping properties of
${\mathfrak L}_\e$ as the scaling parameter $\e$ tends to $0$.
Unfortunately, if we work with the function spaces ${\cal D}^{\ell,
\alpha}_{\e, \delta} ({\R}\times N)$, the norm of the inverse of the
above defined operator blows up as $\e$ tends to $0$.

To get a better understanding about the underlying phenomena, we
digress slightly and study the spectrum of the self-adjoint operator
${\mathfrak L}_\e$ acting on functions which are defined on
$[-1,1]\times N$, have $0$ Dirichlet boundary conditions on $\{\pm
1\}\times N$ and $0$ Neumann boundary conditions on $[-1, 1]\times
\del N$ if $\del N$ is not empty. Recall that there is an
eigenfunction decomposition of the functions $w$ defined on $[-1,1]
\times N$ as
\[
w(t,y) = \sum_{j=0}^\infty w_j(t) \, \phi_j(y) ,
\]
where $(\phi_j, \lambda_j)$ are the eigendata of
$\Delta_h$. This splitting induces a splitting of the operator
${\mathfrak L}_\e$ into a sequence of second order ordinary
differential operators
\[
L_{\e, j}  : =  - \e^2 \, \del_t^2 + \e^2 \, \lambda_j +
\frac{1}{2}\, W'(u_\star (\cdot /\e)) ,
\]
acting on functions defined on $[-1, 1]$. The spectrum of ${\mathfrak
L}_\e$ is then the union of the spectra of the operators $L_{\e,j}$.
Namely
\[
\mbox{Spec} \, ({\mathfrak L}_\e ) = \bigcup_{j \in {\N}} \,
\mbox{Spec} \, (L_{\e, j}) .
\]
Furthermore, the spectrum of $L_{\e,j}$ is equal to the spectrum of
$L_{\e,0}$ shifted by $\e^2 \, \lambda_j$.  Now, the spectrum of
$L_{\e, 0}$ is given by
\[
\mbox{Spec} \, L_{\e, 0} : =  \{ \mu_{0,\e}  < \mu_{1, \e} <  \mu_{2,
\e} \ldots \}.
\]
All the eigenvalues are simple since the operator is a second order
ordinary differential operator. Moreover, arguing as in the proof of
Lemma~\ref{le:fir}, one can show that $\mu_{0, \e} >0$. The existence
of the function $w_\star$ also implies easily that $\mu_{0, \e}$
decays exponentially fast to $0$ as $\e$ tends to $0$. Furthermore,
if  $w_{0, \e}$ denote the eigenfunctions associated to $\mu_{0,\e}$,
which are normalized so that $w_{0, \e}(0)=1$, the sequence of
rescaled functions $(w_{0, \e} (\e \, \cdot))_{\e >0}$ converge on
compact to $w_\star$, as $\e$ tends to $0$. Concerning the second
eigenvalue $\mu_{1, \e}$, one can show that
\[
\lim_{\e\rightarrow 0} \mu_{1, \e} > 0 .
\]
To prove this fact, one can argue by contradiction as in the proof of
Proposition~\ref{pr:8.1}. Since these results are not needed in our
analysis, we leave the proofs to the reader.

On a heuristic level, the above study shows that, as $\e$ tends to
$0$, the number of small eigenvalues of ${\mathfrak L}_\e$ which are
close to $0$ tends to $+\infty$. Moreover the corresponding
eigenfunctions are of the form  $w_{0,\e} \, \phi_j$, where $w_{0,
\e}(\e \cdot)$ tends to $w_\star$ as $\e$ tends to $0$. Since we want
to work with an operator whose inverse is uniformly bounded as $\e$
tends to $0$, we need to replace the spaces we work with by much
smaller function spaces. It is standard to work orthogonally to all
the eigenfunctions which correspond to small eigenvalues to recover
an operator whose inverse is uniformly bounded. Since the
eigenfunctions corresponding to small eigenvalues tend to  $w_\star
(\cdot /\e)  \, \phi_j$,  as $\e$ tends to $0$, it should be  enough
to work orthogonally to all the functions of the form  $w_\star
(\cdot /\e)  \, \phi_j$,  in order to recover an operator whose
inverse is uniformly bounded.

The above should be enough to justify why we replace the spaces
${\cal D}^{\ell, \alpha}_{\delta ,\e} ({\R}\times N)$ by the much
smaller function spaces which we now describe.
\begin{definition}
Given, $\ell \in {\N}$, $\alpha \in (0,1)$, $\delta : =  (\delta_-,
\delta_+) \in {\R}^2$ and $\e \in (0,1)$ such that $\delta_\pm <
\gamma_{\pm}$, we define the weighted space ${\cal E}^{\ell,
\alpha}_{\delta,\e} ({\R}\times N)$ to be the space of functions $u
\in {\cal C}^{\ell, \alpha}_{\delta, \e} ({\R}\times N)$, which
satisfy
\begin{equation}
\forall \, y \in N, \qquad \int_{\R} u (t,y) \, w_\star (t/\e )  \, dt =0 .
\label{eq:9.2}
\end{equation}
This space is endowed with the induced norm.
\end{definition}
Observe that we now impose infinitely many constraints on the
functions $u$, so ${\cal E}^{\ell, \alpha}_{\delta,\e} ({\R}\times
N)$ has infinite codimension in ${\cal C}^{\ell, \alpha}_{\delta,\e}
({\R}\times N)$. The restriction $\delta_\pm < \gamma_\pm$ is needed
to ensure the convergence of the integrals in (\ref{eq:9.2}).

In the case where $N$ has a boundary, we define $\left[ {\cal
E}_{\delta,\e}^{\ell, \alpha}({\R}\times N)\right]_{0}$,
to be the subspace of functions of  ${\cal E}_{\delta, \e}^{\ell,
\alpha}({\R}\times N)$ which have $0$ Neumann boundary condition on
${\R} \times \del N$. Again, we keep the same notation for this space
whether $N$ has a boundary or not.

\subsection{Mapping properties}

It will be interesting to take into account the parameter $\e$ in
order to relax slightly the assumptions on the metric $h$. Indeed, we
will still consider that $h$ satisfies (\ref{eq:8.22}) but instead of
(\ref{eq:8.222}) we now assume that $\ell \geq 2$ is fixed and that
$h$ satisfies the weaker condition
\begin{equation}
||h||_{{\cal C}^{\ell, \alpha}_\e (N)}\leq \frac{1}{\Lambda},
\label{eq:8.zzz}
\end{equation}
where by definition
\begin{equation}
||h||_{{\cal C}^{\ell, \alpha}_\e (N)} : = \| h_{ij} \|_{L^\infty
(N)} + \sum_{j=1}^{\ell-1} \e^{j-1} \, \| \nabla^j h_{ij}
\|_{L^\infty (N)} + \e^{\ell - 2 + \alpha}  \, \| \nabla^{\ell-1}
h_{ij} \|_{{\cal C}^{0, \alpha} (N)}.
\label{eq:8.2222}
\end{equation}

\begin{remark}
At first glance, this condition does not seem very natural since,
paralleling the definition of the spaces ${\cal C}^{\ell, \alpha}_\e
({\R}\times N)$,  it would have been more natural to require that
\begin{equation}
|h|_{{\cal C}^{\ell, \alpha}_\e (N)} \leq \frac{1}{\Lambda},
\label{eq:8.zzzw}
\end{equation}
where
\[
|h|_{{\cal C}^{\ell, \alpha}_\e (N)} : =  \sum_{j=0}^{\ell-1} \e^j \,
\| \nabla^j h_{ij} \|_{L^\infty (N)} + \e^{\ell -1+\alpha}  \, \|
\nabla^{\ell-1} h_{ij} \|_{{\cal C}^{0, \alpha} (N)}.
\]
However, as we will see in the proof of the next Proposition, we will
need a control on the ${\cal C}^1$ norm of $h$ which is slightly
better then the one which is provided by (\ref{eq:8.zzzw}).

Indeed, if $\Lambda$ is assumed to be fixed, $\e \, \nabla h_{ij}$
converges uniformly to $0$ as $\e$ tends to $0$ when $h$ is assumed
to satisfy (\ref{eq:8.zzz}), while this fact is not guaranteed by
(\ref{eq:8.zzzw}). Now the fact that $\e \, \nabla h_{ij}$ converges
uniformly to $0$ as $\e$ tends to $0$ will be used in the proof of
the next Proposition and this will allow us to use the result of
Corollary~\ref{le:6.3}. It would be possible to use (\ref{eq:8.zzzw})
instead of (\ref{eq:8.zzz}) and this would simplify some of the
forthcoming statements. However, we would have to pay a price, namely
we would have to prove a result similar to Corollary~\ref{le:6.3}
when $\Delta$ is replaced by a more general elliptic operator on
${\R}^n$.
\label{re:ddd}
\end{remark}

Observe that, granted the above definition the operator
\begin{equation}
{\mathfrak L}_\e : [{\cal E}^{\ell, \alpha}_{\delta,\e} ({\R}\times
N)]_{0}  \longrightarrow {\cal E}^{\ell-2, \alpha}_{\delta, \e}
({\R}\times N),
\label{eq:9.3}
\end{equation}
is well defined. Indeed, given $w \in {\cal E}^{\ell,
\alpha}_{\delta, \e}({\R}\times N)$ we have
\begin{equation}
\begin{array}{rlllll}
\ds \int_{\R} {\mathfrak L}_\e w \, w_\star(\cdot /\e) \, dt & = &
\ds \int_{\R} \left( - \e^2 (\del_t^2 + \Delta_h) + \frac{1}{2} W''(
w_\star(\cdot /\e))\right) w  \, w_\star(\cdot /\e) \, dt \\[3mm]
& = & \ds \int_{\R} \left( - \e^2\del_t^2+ \frac{1}{2} W''(
w_\star(\cdot /\e)) \right)  w \, w_\star(\cdot /\e) \, dt \\[3mm]
&   & \ds \qquad - \e^2 \,  \Delta_h \left( \int_{\R} w \,
w_\star(\cdot /\e) \, dt \right) \\[3mm]
& = & \ds \int_{\R} \left( - \e^2\del_t^2  + \frac{1}{2}
W''(w_\star(\cdot /\e))\right) w_\star(\cdot /\e) \, w \, dt \\[3mm]
& = &  0 .
\end{array}
\label{eq:XXX}
\end{equation}
Hence, ${\mathfrak L}_\e \, w \in {\cal E}^{\ell-2,
\alpha}_\e({\R}\times N)$. Recall that we have assumed that the
metric $h$ satisfies (\ref{eq:8.22}) and (\ref{eq:8.zzz}). This
implies that ${\mathfrak L}_\e$ is uniformly bounded as $\e$ varies
in $(0, 1)$. Applying the results of \S 8, we see that ${\mathfrak
L}_\e$, defined in (\ref{eq:9.3}) is an isomorphism, provided
$\delta_\pm \in (-\gamma_\pm, \gamma_\pm)$. In the next result, we
show that the inverse of ${\mathfrak L}_\e$ is also uniformly bounded
as $\e$ varies in $(0,1)$.
\begin{proposition}
Assume that $h$ satisfies (\ref{eq:8.22}),(\ref{eq:8.zzz}) and
further assume that $\delta_\pm \in (-\gamma_\pm, 0]$. There exists a
constant $c >0$ (only depending on $\Lambda$), such that, for all $\e
\in (0,1)$ and for all $w \in [{\cal E}^{\ell, \alpha}_{\delta, \e}
({\R}\times N)]_{0}$, we have
\[
||w||_{{\cal C}^{\ell, \alpha}_{\delta, \e} ({\R}\times N)} \leq c \,
|| {\mathfrak L}_\e \, w||_{{\cal C}^{\ell-2, \alpha}_{\delta, \e}
({\R}\times N)}.
\]
\label{pr:9.1}
\end{proposition}
The restriction $\delta_\pm \leq 0$ will be needed to apply the
result of Corollary~\ref{le:6.3}.

{\bf Proof~:} We claim that, there exists a constant $c >0$, such
that if $w$ and ${\mathfrak L}_\e \, w$ are functions bounded by
$\varphi_\delta (\cdot/ \e)$ (and $w$ has $0$ Neumann boundary
condition on ${\R}\times \del N$ if $N$ has a boundary) and if $w$
satisfies (\ref{eq:9.2}), then
\[
||\varphi_{-\delta} (\cdot/ \e) \, w||_{L^\infty ({\R}\times N)} \leq
c \, ||\varphi_{-\delta} (\cdot/ \e) \, {\mathfrak L}_\e \,
w||_{L^\infty ({\R}\times N)}.
\]

As in (\ref{eq:hhh}) in the proof of Proposition~\ref{pr:8.1}, the
proof of the claim relies on the observation that there exists
$t_0 >0$ and $c >0$ such that, if $w$ and ${\mathfrak L}_\e \, w$ are
functions bounded by $\varphi_\delta (\cdot/ \e)$ (and $w$ has $0$
Neumann boundary condition on ${\R}\times \del N$ if $N$ has a
boundary), then
\begin{equation}
|| \varphi_{-\delta} (\cdot/ \e) \, w ||_{L^\infty ({\R} \times N)}
\leq c \left( || w ||_{L^\infty ([-t_0 \e, t_0 \e] \times N)} +  ||
\varphi_{-\delta} (\cdot/ \e) \, {\mathfrak L}_\e \, w ||_{L^\infty
({\R} \times N)} \right).
\label{eq:aure}
\end{equation}
In order to prove the claim, we argue by contradiction. If the claim
where not true, there would exist a sequence $\e_k$ tending to $0$, a
sequence of metrics $h_k$ satisfying (\ref{eq:8.22}) and
(\ref{eq:8.zzz}), a sequence of functions $w_k$ satisfying
(\ref{eq:9.2}) such that
\[
\lim_{k \rightarrow \infty} || \varphi_{-\delta} (\cdot/ \e_k) \,
{\mathfrak L}_{\e_k} \, w_k ||_{L^\infty ({\R}\times N)} = 0 ,
\]
and
\[
|| \varphi_{-\delta} (\cdot/ \e_k) w_k ||_{L^\infty ({\R}\times N)} = 1 .
\]
Moreover, the functions $w_k$ have $0$ Neumann boundary condition on
${\R}\times \del N$ if $N$ has a boundary.

Observe that (\ref{eq:aure}) implies that, for each $k \in {\mathbb N}$
\[
1 \leq  c \left( || \varphi_{-\delta} (\cdot/ \e_k) \, w_k
||_{L^\infty ([- t_0\e_k, t_0\e_k]\times N) } + || \varphi_{-\delta}
(\cdot/ \e_k) \, {\mathfrak L}_{\e_k} \, w_k ||_{L^\infty ({\R}\times
N)} \right).
\]
Furthermore, since the functions $y \longrightarrow \sup_{[-t_0 \e_k,
t_0 \e_k]} \varphi_{-\delta} (t/ \e_k) |w_k (t, y)|$ are continuous,
for each $k \geq 0$, one can choose a point $y_k \in N$ such that
\begin{equation}
1 \leq c \left( || \varphi_{-\delta} (\cdot/ \e_k) \, w_k (\cdot,
y_k) ||_{L^\infty ([- t_0 \e_k, t_0 \e_k])} +  ||\varphi_{-\delta}
(\cdot/ \e_k) \, {\mathfrak L}_{\e_k} \, w_k ||_{L^\infty({\R}\times
N)} \right).
\label{eq:9.4}
\end{equation}
The hypersurface $N$ being compact, we can assume without loss of
generality that the sequence $y_k$ converges in $N$.

We define, for all $s \in {\R}$ and for all $z \in T_{N}(y_k)$ close
enough to the origin
\[
\tilde w_k (s,z) : = w_k (\e_k \, s, \exp_{y_k} (\e_k \, z)).
\]
We identify $T_{N}(y_k)$ with ${\R}^n$. The precise set of  $z \in
{\R}^n$ for which $\tilde w_k$ is well defined depends on the
injectivity radius of $N$. In the case where $N$ has a boundary, this
set also depends on the distance from $y_k$ to $\del N$,
\[
\rho_k := d(y_k, \del N) /\e_k ,
\]
where $d$ denotes the geodesic distance on $N$ for the metric $h_k$.
We now distinguish a few cases according to the behavior of the
sequence $(\rho_k)_k$, which can be assumed to converge to
$\rho_\infty \in [0, \infty]$.

\noindent
{\bf Case 1.} Assume that either $N$ has no boundary or $N$ has a
boundary and  $\rho_\infty = +\infty$. Since both $w_k$ and
${\mathfrak L}_{\e_k} \, w_k $ are uniformly bounded on compact sets,
we deduce from elliptic estimates that the sequence of functions
$\tilde w_k$ is bounded in ${\cal C}^{1, \alpha}$ topology on any
compact of ${\R}\times {\R}^n$. Now, Ascoli's Theorem together with a
standard diagonal argument implies that, up to a subsequence and for
some $\alpha' < \alpha$ the sequence of functions $\tilde w_k$
converge (in ${\cal C}^{1, \alpha'}$ topology) to $\tilde w_\infty$.
Furthermore, $w_\infty$ is bounded by $\varphi_\delta$. Passing to
the limit in the equation satisfied by $\tilde w_k$, we claim that
the function $\tilde w_\infty$ is a weak solution of
\[
\left( \del_t^2 + \Delta - \frac{1}{2} \, W'' (u_\star) \right) \,
\tilde w_\infty = 0 ,
\]
in ${\R}\times {\R}^n$. This is precisely at this point that we use
the full strength of (\ref{eq:8.zzz}) which guarantees that the
${\cal C}^1$ norm of the coefficients of the metrics $h_k$ are
uniformly bounded, hence the sequence of dilated metrics converges to
the {\em flat} metric. While the weaker more natural condition given
in Remark~\ref{re:ddd} would only have ensured that the sequence of
dilated metrics converges to {\em some} metric on ${\R}^n$. Indeed,
observe that, in order to define the functions $\tilde w_k$, we have
chosen normal coordinates $\e_k \, z$ to parameterize $N$ in some
neighborhood of the point $y_k$. In these coordinates the metric
$h_k$ can be written as
\[
h_k = \e^2_k \, \tilde h_{k,jj'} (z) \, dz_j \otimes dz_{j'},
\]
The fact that $\e_k \, z$ are normal coordinates implies that
\[
\tilde h_{k,jj'} (0) = \delta_{jj'},
\]
and (\ref{eq:8.zzz}) translates into
\[
||\nabla^m \tilde h_{k,jj'}||_{L^\infty} \leq c \, \e_k,
\]
for $m=1, \ldots, \ell$. Here all norms are evaluated on a fixed
compact, $c>0$ only depends on $\Lambda$ and the partial derivatives
are taken with respect to $z_i$. This in turn implies that the
Laplace-Beltrami operator $\e^2_k \, \Delta_{h_k}$ can be expanded as
\[
\e^2_k \, \Delta_{h_k} = \Delta + {\cal O} (\e_k \, |z|\,
\del_{z_j}\, \del_{z_{j'}}) + {\cal O} (\e_k \, \del_{z_j} ) ,
\]
where $\Delta$ is the ordinary Laplace-Beltrami operator on ${\R}^n$.
This explains why the equation satisfied by $\tilde w_\infty$
involves the operator $\Delta$.

\noindent
{\bf Case 2.} Assume that $N$ has a boundary and that  $\rho_\infty <
+\infty$. Applying the above argument, we can assume that, up to a
subsequence, the sequence of functions $\tilde w_k$ converges to
$\tilde w_\infty$. Furthermore, $\tilde w_\infty$ is a  weak solution
of
\[
\left( \del_t^2 + \Delta - \frac{1}{2} \, W'' (u_\star) \right) \,
\tilde w_\infty = 0 ,
\]
in ${\R}\times {\R}^n_+$ which is bounded by $\varphi_\delta$. Here,
we have set ${\R}^n_+ : =  \{ (z_1, \ldots, z_n) \, : \, z_n >0 \}$
and we have assumed that the tangent space at a point of $\del N$ is
identified with ${\R}^n_+$. In addition, $w$ has $0$ Neumann boundary
data and we  can then extend $\tilde w_\infty$ to ${\R}\times {\R}^n$
by reflection across the hyperplane $z_n=0$ to reduce to Case 1.

In either case, Lebesgue's convergence Theorem implies that $\tilde
w_\infty$ satisfies
\begin{equation}
\int_{{\R}\times {\R}^n} \tilde w_\infty \, w_\star \, dt\, dz =0 .
\label{eq:9.5}
\end{equation}
We give the proof of this fact when $\rho_\infty = +\infty$, obvious
modifications are needed to handle the case where $\rho_\infty <
+\infty$. Given a function $\eta : {\R}^n \longrightarrow {\R}$, with
compact support we use the fact that, for all $z$ in ${\R}^n$ and all
$k$ large enough
\[
\int_{\R} \tilde w_k (s,z) \, w_\star (s) \, ds =0.
\]
Provided $k$ is chosen large enough so that the domain of definition
of $z \rightarrow \tilde w_k(s,z)$ includes the support of $\eta$, we
multiply this identity by $\eta$ and integrate over $z$ to get
\[
\int_{\R \times N} \eta (z) \, \tilde w_k (s,z) \, w_\star (s) \, ds
\, da_{h_k} = 0 .
\]
Passing to the limit as $k$ tends to $\infty$, we conclude that
\[
\int_{\R \times N} \eta (z) \, \tilde w_\infty (s,z) \, w_\star (s)
\, ds \, dz = 0 .
\]
This identity being valid for all $\eta$, we finally obtain, for all
$z \in {\R}^n$
\[
\int_{\R} \tilde w_\infty (s,z) \, w_\star (s) \, ds =0 .
\]
By construction, the function $\tilde w_\infty$ is bounded by
$\varphi_\delta$ and since we have assumed that $\delta_\pm \leq 0$,
this implies that $\tilde w_\infty$ is bounded. We now apply the
result of Corollary~\ref{le:6.3} which shows that $\tilde w_\infty$
only depends on $t$. Since all bounded solutions of $L_0 \, w=0$ are
collinear to $w_\star$, we can write $\tilde w_\infty (s, z) =  c \,
w_\star (s)$ for some constant $c\in {\R}$, and (\ref{eq:9.5})
implies that $c =0$. Therefore, we have shown that $\tilde w_\infty
\equiv 0$.

Finally, (\ref{eq:9.4}) implies that, for each $k \in {\mathbb N}$
\[
1 \leq c \, \left( || \varphi_{-\delta}  \, \tilde w_k (\cdot,
0)||_{L^\infty( [- t_0 , t_0 ] \times N)} +  || \varphi_{-\delta}
(\cdot/\e_k) \, {\mathfrak L}_{\e_k} \, w_k ||_{L^\infty ({\R}\times
N)} \right),
\]
passing to the limit as $k$ tends to $\infty$ in this inequality, we
conclude that
\[
1 \leq c \, || \varphi_{-\delta} \, \tilde w_\infty (\cdot, 0)
||_{L^\infty( [- t_0 , t_0 ] \times N)} ,
\]
which implies that $\tilde w_\infty$ is not identically equal to $0$.
Since we have reached a contradiction, the proof of the claim is
complete.

The claim being proved the result of the Proposition follows
immediately from Schauder's estimates which are applied in geodesic
balls of radius $\e$ in ${\R}\times N$. \hfill $\Box$

Obviously, this result could have been obtained directly without any
reference to the result of Proposition~\ref{pr:8.1}. However, we feel
that the decomposition of the proof into two different steps sheds
light on the choice of the function spaces and in particular explains
where the conditions (\ref{eq:8.2}) and (\ref{eq:9.2}) enter into
play.

\section{Fermi coordinates}\setcounter{equation}{0}

Assume that $N \subset M$ is an admissible hypersurface. To begin
with, let us assume that $N$ has no boundary. We will explain in the
last paragraph of this section the modifications which are needed to
handle the case where the boundary of $N$ is not empty.

For any $p \in M$ we define $d_{N} (p)$ to be the signed geodesic
distance from $p$ to $N$. This means that $d_{N} (p)$ is the geodesic
distance from $p$ to $N$ if $p \in M^+ (N)$ and that $d_N (p)$ is
equal to minus the geodesic distance from $p$ to $N$ if $p \in M^-
(N)$. We define
\begin{equation}
V_\tau (N) : = \{p \in M \quad : \quad d_{N}(p) \in (-\tau, \tau)\} .
\label{eq:12.1}
\end{equation}
It is well known that for all $\tau$ small enough, the set $V_\tau
(N)$ is a tubular neighborhood of $N$ in $M$, provided $N$ is at
least ${\cal C}^2$.  In order to measure the regularity of a function
$u$ which is defined in $V_{\tau} (N)$, we define for $\ell \geq 0$,
$\alpha \in (0,1)$ and $\e \in (0,1)$
\[
||u||_{{\cal C}^{\ell, \alpha}_{\e} (V_\tau (N))} : = \sum_{j=0}^\ell
\e^j \, \|\nabla^j u \|_{L^\infty (V_\tau (N))}+ \e^{\ell+\alpha} \,
\sup_{p\neq q \, \in \, V_\tau (N)} \frac{|\nabla^\ell u(p)-
\nabla^\ell u(q)|}{d(p,q)^\alpha} ,
\]
where $d$ is the geodesic distance in $M$.

In the previous sections, we have asked that the metrics $h$ we
consider on the manifold $N$ satisfy (\ref{eq:8.22}) and
(\ref{eq:8.222}) or (\ref{eq:8.zzz}). This was needed to guarantee
that the results we have obtained are independent of the choice of
$h$, but might depend on $\Lambda$. We now restrict our attention to
admissible hypersurfaces $N$ which are embedded in $M$ and which are
close to some reference smooth admissible hypersurface $N_0$. To make
things precise, we consider a smooth admissible hypersurfaces $N_0$
and assume that all the admissible hypersurfaces we consider can be
written as a geodesic normal graph over $N_0$ for some function
$\psi_N  \in {\cal C}^{\ell, \alpha} (N_0)$.
\begin{definition}
The ${\cal C}^{\ell, \alpha}_\e$ norm of the hypersurface $N$, with
respect to the hypersurface $N_0$, is defined to be
\[
||N||_{{\cal C}^{\ell, \alpha}_\e (N_0)}: =  ||\psi_N||_{{\cal
C}^{\ell, \alpha}_\e (N_0)},
\]
where
\begin{equation}
\begin{array}{rllll}
||\psi||_{{\cal C}^{\ell, \alpha}_\e (N_0)} & := & \displaystyle \|
\psi \|_{L^\infty (N_0)} + \| \nabla \psi \|_{L^\infty (N_0)} \\[3mm]
& + & \displaystyle \sum_{j=0}^{\ell-2} \e^{j}  \| \nabla^{j+2} \psi
\|_{L^\infty (N_0)} + \e^{\ell - 2 + \alpha} \| \nabla^{\ell} \psi
\|_{{\cal C}^{0, \alpha} (N_0)} ,
\end{array}
\label{eq:10.10}
\end{equation}
and where $\psi_N$ is the function whose normal geodesic graph over
$N_0$ is $N$.
\label{def:nnn}
\end{definition}
Since the hypersurface $N$ is assumed to be embedded in $M$, the
metric $h$ we consider on $N$ is just the metric induced by the
metric $g$ of $M$. Now, given $\Lambda >0$ and $\ell \geq 2$, one can
find constants $\Lambda', \Lambda'' >0$ such that~:
\begin{itemize}
\item[(i)] If the ${\cal C}^{2} (N_0)$ norm of the function $\psi_N$,
whose normal geodesic graph is $N$, is bounded by $\Lambda'$, then
$N$ is embedded.
\item[(ii)] If the ${\cal C}^{\ell, \alpha}_\e(N_0)$ norm of $\psi_N$
is bounded by $\Lambda''$, then conditions (\ref{eq:8.22}) and
(\ref{eq:8.zzz}) are fulfilled.
\end{itemize}

\begin{remark}
As in (\ref{eq:8.2222}) the definition of the weighted norm of a
hypersurface is not the natural one, namely
\[
|N|_{{\cal C}^{\ell, \alpha}_\e (N_0)} : =  |\psi_N|_{{\cal C}^{\ell,
\alpha}_\e (N_0)},
\]
where
\begin{equation}
|\psi|_{{\cal C}^{\ell, \alpha}_\e (N_0)} :=  \displaystyle
\sum_{j=0}^{\ell} \e^{j} \,  \| \nabla^j \psi \|_{L^\infty (N_0)} +
\e^{\ell + \alpha} \, \| \nabla^{\ell} \psi \|_{{\cal C}^{0, \alpha}
(N_0)} .
\label{eq:10.1010}
\end{equation}
In order to justify this choice, observe that, as already mentioned,
Definition~\ref{def:nnn} is consistent with the definition of the
norm of the induced metric which is used in \S 9.2, while the above
definition is consistent with the definition used in
Remark~\ref{re:ddd}.
\label{re:dddd}
\end{remark}

In order to define Fermi coordinates in some neighborhood of $N$ we
set for all $t \in [-\tau, \tau]$ and all $y \in N$
\[
Z_N (t,y) : = \exp_y (t \, \nu_y),
\]
where $\nu_y$ is the normal to $N$ at the point $y$, which is assumed
to point towards $M^+(N)$. The parameter $\tau$ is chosen small
enough so that $Z_{N}$ is a diffeomorphism onto its image.
\begin{definition}
The coordinates $(t,y)$, defined in $V_\tau (N)$, are called Fermi
coordinates relative to the hypersurface $N$.
\end{definition}
We will need the following classical result which states that the
distance function to $N$ is well defined and as smooth as $N$ in some
tubular neighborhood of $N$ whose width is bounded from below. The
proof of this result can be found in \cite{Kra-Par} or in \cite{gt}.
\begin{lemma}
Assume that $\ell \geq 2$, $\Lambda' >0$ is fixed small enough and
$\Lambda'' >0$ is fixed. Then, there exists $\tau_0 >0$ only
depending on $\Lambda'$ such that if $N$ is a hypersurface whose
${\cal C}^{2}(N_0)$ norm is bounded by $\Lambda'$ and whose ${\cal
C}^{\ell, \alpha}_\e (N_0)$ norm is bounded by $\Lambda''$, then
$Z_N$ is a ${\cal C}^{\ell-1, \alpha}$ diffeomorphism from $(-\tau_0,
\tau_0)\times N$ onto its image. Moreover the norm of $\nabla d_N$ in
${{\cal C}^{\ell-1, \alpha}_\e (V_{\tau_0} (N))}$ is bounded by a
constant which only depends on $\Lambda''$.
\label{le:11.1}
\end{lemma}
{\bf Proof :} The existence of the constant $\tau_0$ follows from the
fact that the hypersurfaces we consider are uniformly bounded in
${\cal C}^{2}$ topology. This implies that the principal curvatures
of $N$ are bounded from above by some constant only depending on
$\Lambda'$ and, in turn, this shows that the size of the tubular
neighborhood over which $Z_N$ is a diffeomorphism is bounded from
below.

The fact that $d_N$ is as regular as $N$ follows from standard
arguments. Indeed, the gradient of $d_N$ at the point of coordinates
$(t,y)$ is given by the tangent vector to the geodesic starting from
$y$ with vector speed $\nu_y$. Hence the gradient of $d_N$ is as
regular as the Gauss map of $N$, which in turn is as regular as
$\nabla \psi_N$. Hence $\nabla d_N$ is bounded in ${\cal C}^{\ell-1,
\alpha}_\e(V_{\tau_0})$ by $\Lambda''$ (we even know that $\nabla^2
d_N$ is bounded in ${\cal C}^{\ell-2, \alpha}_\e(V_{\tau_0})$ by some
constant only depending on $\Lambda''$, which is a slightly stronger
statement). \hfill $\Box$

Let us recall that
\begin{equation}
|\nabla \, d_N |_g =1,
\label{eq:11.1}
\end{equation}
a.e. in $M$.

>From now on, we assume that $\tau_0$ is chosen as in
Lemma~\ref{le:11.1}. For all $t \in [-\tau_0, \tau_0]$, we define the
hypersurface $N_t$ to be the hypersurface which is parallel to $N$,
at distance $t$. This means that $N_t$ is the normal geodesic graph
over $N$ for the constant function $\psi(y)\equiv t$. This definition
being understood, we recall the expression of the Laplace-Beltrami
operator in Fermi coordinates.
\begin{lemma}
Denote by $(t,y)$ the Fermi coordinates relative to $N$, which are
defined in $V_{\tau_0}(N)$. Then, at any point of $Z_N(t, y) \in N_t$
we have
\[
\Delta_g = \del_t^2 + \Delta_{h_t} - n \, H_{N_t} \,\del_t .
\]
where $h_t$ is the metric induced by $g$ on $N_t$ and where
$H_{N_t}(t,y)$ denotes the mean curvature of the hypersurface $N_t$
at the point $Z_N(t,y)$.
\label{le:11.2}
\end{lemma}
{\bf Proof :} We provide the proof of the formula when $t=0$. The
general case can be treated similarly, up to notation changes. Let
$e_1, \ldots, e_n$ be an orthonormal frame field on $N$ and $\nu$ be
the normal vector field.

The Laplace-Beltrami operator on $M$ is defined by
\[
\Delta_g =  \sum_{i=1}^n \left( e_i \, e_i  - D_{e_i} e_i \right)  +
\nu  \,  \nu - D_\nu \nu ,
\]
where $D$ is the Levi-Civita connection on $M$. Let $D^N$ denote the
Levi-Civita connection on $N$, by construction, we have
\[
D_{e_i} e_i = D^N_{e_i} e_i + g( D_{e_i} e_i , \nu)  \, \nu .
\]
Therefore
\[
\Delta_g =  \sum_{i=1}^n \left( e_i \, e_i  - D^N_{e_i} e_i \right)
+ \sum_{i=1}^n  g(e_i , D_{e_i} \nu)  \, \nu  +  \nu  \,  \nu - D_\nu
\nu .
\]
By definition $\nu\, \nu = \del_t^2$ and $\nu = \del_t$. Furthermore
   $D_\nu \nu=0$ and
\[
\sum_{i=1}^n  g(e_i , D_{e_i} \nu) = - n \, H_N ,
\]
where $H_N$ is the mean curvature of $N$. Hence we conclude that
\[
\Delta_g =   \del_t^2 + \Delta_h - n \, H_N \, \del_t .
\]
Which is the desired expression. \hfill $\Box$

Applying the expansion given in the previous Lemma to the function
$d_N (t,y)=t$, we obtain the important formula
\begin{equation}
\Delta_g \, d_N (t,y)  = - n \, H_{N_t}(t,y).
\label{eq:11.2}
\end{equation}
When $t=0$ we recover the well known fact that the Laplacian of the
distance function to the hypersurface $N$, computed at a point of
$N$, is equal to (minus) the sum of the principal curvatures of the
hypersurface at this point \cite{gt}.

We will also need the expansion of the Laplace-Beltrami operator
$\Delta_g$ in Fermi coordinates. To do so, we need slightly better
control on the regularity of the hypersurfaces $N$ which should be at
least  bounded in ${\cal C}^{3, \alpha}_\e(N_0)$.
\begin{proposition}
Assume that $\ell\geq 2$, $\Lambda' >0$ is fixed small enough and
$\Lambda'' >0$ is fixed. If $N$ is a hypersurface whose  ${\cal
C}^{2}(N_0)$ norm is bounded by $\Lambda'$ and whose ${\cal
C}^{\ell+1, \alpha}_\e (N_0)$ norm is bounded by $\Lambda''$, then
there exists a second order operator in $\del_{y_j}$
\[
L_2 : = t \, {\cal O}_{{\cal C}^{\ell-1, \alpha}_\e} (\Lambda'') \, \nabla^2,
\]
without any first or zero-th order term, and a first order operator
in $\del_t$ and $\del_{y_j}$
\[
L_1 := {\cal O}_{{\cal C}^{\ell-2, \alpha}_\e}(\Lambda'') \, \nabla +
t \, \e^{-1} \, {\cal O}_{{\cal C}^{\ell-2, \alpha}_\e} (\Lambda'')
\, \nabla,
\]
without zero-th order term,  such that
\[
\Delta_g =  \del^2_t + \Delta_{h}  + L_1 + L_2 ,
\]
where $h$ is the metric on $N$, induced by the metric $g$ on $M$. The
notation ${\cal O}_{{\cal C}^{k, \alpha}_\e} (\Lambda'')$ refers to
the fact that the coefficients of the operators are bounded in
${\cal C}^{k, \alpha}_\e (V_{\tau_0}(N))$ by a constant which only
depends on $\Lambda''$.
\label{co:11.2}
\end{proposition}
{\bf Proof~:} This result follows from the expansion of the metric in
Fermi coordinates $(t,y)$ relative to $N$. Indeed, if we write the
metric $g$ as
\[
g = g_{tt}\, dt\otimes dt+ \sum_{j=1}^n g_{tj}\, dt \otimes dy_j +
\sum_{j,j'=1}^n  g_{jj'} \, dy_j \otimes dy_{j'},
\]
the coefficients of the metric can be expanded as
\[
g_{tt} = 1 , \qquad  \qquad  g_{tj} = 0 \qquad \mbox{and} \qquad
g_{jj'} = h_{jj'} + t \, {\cal O}_{{\cal C}^{\ell-1, \alpha}_\e}
(\Lambda'') ,
\]
where $h_{jj'}$ are the coefficients of the induced metric on $N$.
This is precisely at this point that we use the definition
(\ref{eq:10.10}). Indeed, the regularity of the coefficients of the
metric $g$ in Fermi coordinates is the same as the regularity of the
gradient of the Gauss map of $N$, which in turn can be estimated by
the second derivatives of the function whose graph is $N$. Hence, the
${\cal C}^{\ell-1, \alpha}_\e (V_{\tau_0}(N))$ norm of the
coefficients of the metric is bounded by the ${\cal C}^{\ell+1,
\alpha}_\e(N_0)$ norm of $N$.

Similar expansions hold for $g^{\alpha\beta}$, the coefficients of
$g^{-1}$. This implies that
\[
\sqrt{\mbox{det } g} =  \sqrt {\mbox{det } h} +  t\, {\cal O}_{{\cal
C}^{\ell-1, \alpha}_\e} (\Lambda'').
\]
With these expansions, we get
\[
\begin{array}{rlllll}
\ds \frac{1}{\sqrt{\mbox{det } g}} \, \del_t \left( \sqrt{\mbox{det }
g} \, g^{tt} \, \del_t \right) & = & \del_t^2 + \left( {\cal
O}_{{\cal C}^{\ell-2, \alpha}_\e} (\Lambda'')  +  t \, \e^{-1}\,
{\cal O}_{{\cal C}^{\ell-2, \alpha}_\e} (\Lambda'') \right) \, \del_t,
\end{array}
\]
and
\[
\begin{array}{rlllll}
\ds \frac{1}{\sqrt{\mbox{det } g}} \, \del_{y_j} \left(
\sqrt{\mbox{det } g} \, g^{jj'} \, \del_{y_{j'}} \right) & = & \ds
\frac{1}{\sqrt{\mbox{det } h}} \, \del_{y_j} \left( \sqrt{\mbox{det }
h} \, h^{jj'} \, \del_{y_{j'}} \right) \\[3mm]
       & + & t \, \e^{-1}\,  {\cal O}_{{\cal C}^{\ell-2, \alpha}_\e}
(\Lambda'') \, \del_{y_j} + t \, {\cal O}_{{\cal C}^{\ell-1,
\alpha}_\e} (\Lambda'') \, \del_{y_j} \, \del_{y_{j'}} .
\end{array}
\]
The result now follows at once. \hfill $\Box$

As promised, we now explain the modifications needed when $\ptl
N\neq\emptyset$.  Observe that in this case the geodesics starting
from a point $y \in N$ with initial velocity $\nu_{N}(y)$ induce a
fibration of a tubular neighborhood of
$N$.  In the case where $N$ has a boundary this property is not true
anymore.  To overcome this problem, we define what we will
call "twisted" Fermi coordinates.  Recall that we have
assumed in the introduction that $M$ is a smooth domain of a
Riemannian manifold $\tilde M$.  An admissible hypersurface $N$ in $M$
can be extended to a hypersurface $\tilde N$ whose ${\cal C}^{\ell,
\alpha}_\e$ norm is bounded by a constant (independent of $\e$) times
the ${\cal C}^{\ell, \alpha}_\e$ norm of $N$.

Observe that, provided $\tau_0$ is chosen small enough, we can define Fermi
coordinates in
\[
   V_{\tau_0} (\tilde N) : = \{ Z_{\tilde N}(t,z) \in \tilde M \, :
\, z \in \tilde N, \quad t \in (-\tau_0, \tau_0)\}.
\]
Let $X_1$ be the vector field $\ptl_{t}$, which is defined in
$V_{\tau_{0}}(\tilde N)$. Reducing $\tau_{0}$ if this is necessary,
we can assume that $\nu_{\ptl M}$, the inner normal to $\ptl M$, does
not coincide with $X_1$ in $V_{\tau_{0}}(\tilde N)$ (recall that $N$
is assumed to meet $\ptl M$ in an orthogonal way).

We set $X_2$ to be the extension of $\nu_{\ptl M}$ which is given by
the unit normal vector fields to the family of parallel hypersurfaces
to $\ptl M$. This vector field is well defined in a fixed tubular
neighborhood of radius $\e_{0}>0$ of $\ptl M$. Finally, we consider a
smooth cut-off function $\eta$ which is identically equal to $1$ in
$(-\infty, \e_{0}/2)$ and equal to $0$ in $(\e_{0},\infty)$, and let
\[
\tilde \chi : = \eta (d_{\ptl M}),
\]
where $d_{\del M}$ is the distance function to $\del M$. These
definitions being understood, we define the vector field
\[
X : =\frac{X_1 - \left< X_1, \tilde \chi  \, X_2 \right>_g\, \tilde \chi
X_2}{1-\left< X_1 , \tilde \chi \, X_2\right>_g^2}.
\]
It is straightforward to check the following properties of $X$~:
\begin{itemize}
\item[(i)] $X$ is tangent to $\ptl M$.
\item[(ii)] $X=X_1$ at a point $p$ whenever $d_{\del M} (p)\geqslant \e_{0}$.
\item[(iii)] $\left<X, X_1 \right>_g \equiv 1$.
\end{itemize}
We now define $\{F_{t}\}_{t}$ to be the flow associated to $X$. It
follows from (i) that $F_t (N)$ is contained in $M$. Also, (i) and
(iii) imply that $(t,y)\rightarrow F_t(y)$ is in fact a
diffeomorphism from $(-\tau_{0},\tau_{0})\times N$ onto
$V_{\tau_{0}}(\tilde N)\cap M$.

The ``twisted'' Fermi coordinates $(t,y)$ relative to $N$ are defined by~:
\[
Y_{N}(t,y):= F_{t}(y).
\]
Observe that (ii) implies that $Y_{N}(t,y) = Z_{N}(t,y)$ away from a
tubular neighborhood of radius $\e_{0}$ around $\ptl M$. Moreover,
(iii) together with the fact that, $N$ being an admissible
hypersurface, $N$ meets $\del M$ orthogonally,imply that in a
neighborhood of radius $\e_{0}$ of $\ptl M$, if $(t,y)$ are the
``twisted" Fermi coordinates of a point $p$ (relative to $N$) and if
$(s, z)$ are the Fermi coordinates of the same point $p$ (relative to
$\tilde N$), we have
\[
s=t , \qquad \mbox{and} \qquad z - y = {\cal O}(t^2).
\]
Using these properties, one checks that all the results we have
obtained in this section do hold when Fermi coordinates are replaced
by ``twisted" Fermi coordinates.

One last comment about the modifications which are needed when $\del
N$ is not empty. In the case where $\del N$ is empty, one can
parameterize any hypersurface $\tilde N$ close enough to $N$ as a
normal geodesic graph, that is, using the flow $Z_{N}$. In turn, when
$\del N$ is not empty, we agree that we will parameterize any
hypersurface $\tilde N$ close enough to $N$ using the flow $Y_N$.
With slight abuse of terminology we will say that $\tilde N$ is a
``normal graph'' over $N$.

\section{The linear problem}\setcounter{equation}{0}

Assume that $(M,g)$ is a smooth $(n+1)$-dimensional Riemannian
manifold, with or without boundary and that $N_0 \subset M$ is a
fixed admissible, nondegenerate minimal hypersurface or a
volume-nondegenerate constant mean curvature hypersurface.

\subsection{Preliminary assumptions and definitions}

From now on, we assume that $\ell \geq 2$, $\Lambda' >0$ small enough
and $\Lambda'' >0$ are fixed. We agree that, from now on, all the
admissible hypersurface $N$ we consider, satisfy~:
\begin{itemize}
\item[($H$)]  \qquad The ${\cal C}^{2} (N_0)$ norm of the function
$\psi_N$, whose normal graph is $N$, is bounded by $\Lambda'$.
\item[($H'_\ell$)] \qquad The ${\cal C}^{\ell, \alpha}_\e(N_0)$ norm
of the function $\psi_N$, whose normal graph is $N$, is bounded by
$\Lambda''$.
\end{itemize}
We assume that $\Lambda'$ and $\Lambda''$ are chosen so that whenever
$N$ satisfies the above conditions, then $N$ is embedded and
(\ref{eq:8.22}) and (\ref{eq:8.zzz}) are fulfilled. Furthermore,
reducing $\Lambda'$ if this is necessary, we can assume that $N
\subset V_{\tau_0/4} (N_0)$ where $\tau_0$ is chosen as in
Lemma~\ref{le:11.1}.

We choose a cutoff function $\chi$ which is identically equal to $1$
in $(-\tau_0 /2, \tau_0/2)$ and identically equal to $0$ outside
$[-\tau_0,\tau_0]$. We define
\begin{equation}
u_{\e, N} : =  \chi (d_{N}) \, \frac{d_{N}}{|d_{N}|} + (1 - \chi
(d_{N})) \, u_\star (d_{N} / \e) .
\label{eq:aprox}
\end{equation}

\subsection{Function spaces}

Paralleling what we have done in \S 9, we define weighted H\"older
spaces adapted to our problem. We first give the definition of the
H\"older spaces which take into account the scaling parameter $\e$.
\begin{definition}
Given $\ell' \in {\N}$, $\alpha \in (0,1)$ and $\e \in (0,1)$, we
define the space ${\cal C}^{\ell', \alpha}_{\e} (M)$ to be the space
of functions which are $\ell'$ times differentiable and whose
$\ell'$-th partial derivatives are H\"older of exponent $\alpha$.
This space is endowed with the norm
\[
||u||_{{\cal C}^{\ell', \alpha}_{\e} (M)} : = \sum_{j=0}^{\ell'}
\e^j \, \|\nabla^j u \|_{L^\infty (M)  }+ \e^{\ell'+\alpha} \,
\sup_{p\neq q \, \in \, M} \frac{|\nabla^{\ell'} u(p)- \nabla^{\ell'}
u(q)|}{d(p,q)^\alpha} ,
\]
where $d$ is the geodesic distance in $M$.
\end{definition}

This definition being understood, we define a projection operator~:
\begin{definition}
Denote by $(t, y)$ the (twisted) Fermi coordinates relative to $N$.
For all $u \in {\cal C}^{\ell, \alpha} (M)$, we define the function
$\Pi_{\e, N} (u)$ by
\begin{equation}
\Pi_{\e, N} (u) (t,y) : = u (t,y) - \displaystyle \frac{\displaystyle
\int_{\R} \chi  (s) \, u (s, y) \, w_\star (s /\e) \, ds
}{\displaystyle \int_{\R} \chi^2 (s) \, w_\star^2 (s /\e) \, ds} \,
\chi (t) \, w_\star (t/\e),
\label{eq:13.999}
\end{equation}
in $V_{\tau_0} (N)$ and $\Pi_{\e, N} (u) = u$ in $M - V_{\tau_0}(N)$.
\end{definition}

Observe that $\Pi_{\e, N}$ is an involution since, by construction
$\Pi_{\e, N}\circ \Pi_{\e, N} = \Pi_{\e, N}$. It is straightforward
to check the~:
\begin{lemma}
Assume that $(H)$ and $(H'_{\ell+1})$ hold for some $\ell \geq 2$.
Then, for all $\ell'=0, \ldots, \ell$
\[
\Pi_{\e, N} :  {\cal C}^{\ell', \alpha}_\e (M) \longrightarrow {\cal
C}^{\ell', \alpha}_\e (M),
\]
is well defined and uniformly bounded for $\e \in (0,1)$.
\label{le:proj}
\end{lemma}
{\bf Proof~:} The fact that this operator is well defined follows
from the fact that Fermi coordinates induce a local ${\cal C}^{\ell,
\alpha}$ diffeomorphism ${\cal D}$ when $N$ is a ${\cal C}^{\ell+1,
\alpha}$ hypersurface.

The fact that $\Pi_{\e, N}$ is bounded uniformly between the weighted
spaces follows at once from the definition of the ${\cal C}^{\ell+1,
\alpha}_\e(N_0)$ norm of $N$ which ensures that the diffeomorphism
${\cal D}$ is uniformly bounded in $L^\infty (V_{\tau_0} (N))$ and
has a differential which is uniformly  bounded in ${\cal
C}^{\ell-1}_\e (M)(V_{\tau_0}(N))$. \hfill $\Box$

We now define weighted function spaces on $M$  which parallel the
spaces ${\cal E}^{\ell, \alpha}_{\delta,\e} ({\R}\times N)$, when the
weights $\delta_\pm$ are chosen to be equal to $0$.
\begin{definition}
Assume that $(H)$ and $(H'_{\ell+1})$ hold for some $\ell \geq 2$.
Given $0 \leq \ell' \leq \ell$, $\alpha \in (0,1)$ and $\e \in
(0,1)$, we define the space ${\cal E}^{\ell', \alpha}_{\e, N} (M)$ to
be the space of functions $u \in {\cal C}^{\ell', \alpha}_{\e}(M)$
which satisfy
\[
\Pi_{\e, N} (u) = u .
\]
This space is endowed with the induced norm.
\end{definition}

As usual, we also define $[{\cal C}^{\ell', \alpha}_\e (M)]_0$ (resp.
$[{\cal E}^{\ell', \alpha}_{\e,N} (M)]_0$) to be the subspace of
functions of ${\cal C}^{\ell', \alpha}_\e (M)$ (resp. ${\cal
E}^{\ell', \alpha}_{\e,N} (M)$) which have $0$ Neumann boundary
condition on $\del M$, if this boundary is not empty.

\subsection{The linear problem}

Fix $\ell \geq 2$ and assume that $N$ is an admissible hypersurface
which satisfies $(H)$ and $(H'_{\ell})$. Let $u_{\e, N}$ be defined
as in (\ref{eq:aprox}), the linearized operator we are interested in
reads
\[
{\mathbb L}_{\e, N} : =  - \e^2 \, \Delta_g +  \frac{1}{2}\, W'' (u_{\e, N}) .
\]
In this section, we would like to construct a right inverse for this
operator. To do so, we first define and study the auxiliary operator
\[
{\mathbb L}_0 : =  - \e^2 \, \Delta_g + \frac{1}{2} \, \Gamma ,
\]
where the potential $\Gamma$ is chosen to interpolate smoothly
between $2 \, \gamma_-^2 = W''(-1)$ in $M^-(N)$ and $2 \, \gamma_+^2
=  W''(1)$ in $M^+(N)$. More precisely, let $\xi$ denote a smooth
cutoff function equal to $1$ on $(1, +\infty)$ and equal to $0$ on
$(-\infty, -1)$, with $\xi \geq0$. If $(t,y)$ are (twisted) Fermi
coordinates relative to $N$, we define
\[
\Gamma (t,y) : = 2 \, \left( (1-\xi (t/\e))\, \gamma_-^2 + \xi (t/\e)
\,  \gamma_+^2\right) ,
\]
in $V_{\tau_0} (N)$ and
\[
\Gamma  : =  2 \, \gamma^2_\pm,
\]
in $M^\pm (N)-V_{\tau_0} (N)$. We have~:
\begin{lemma}
Fix $\ell \geq 2$. There exists $\e_0 >0$ only depending on
$\Lambda'$ and $\Lambda''$ such that, for all $\e \in (0, \e_0)$, the
operator
\[
{\mathbb L}_0 : {\cal C}^{\ell-2, \alpha}_\e (M)  \longrightarrow
[{\cal C}^{\ell, \alpha}_\e (M)]_0,
\]
is an isomorphism the norm of whose inverse is bounded by some
constant which depends on $\Lambda', \Lambda''$ but does not depend
on $\e$ nor on $N$ satisfying $(H)$ and $(H'_{\ell})$.
\label{le:fgfg}
\end{lemma}
{\bf Proof~:} The fact that ${\mathbb L}_0$ is an isomorphism,
clearly follows from the fact that the potential $\Gamma$  is bounded
from below by a positive constant independent of $\e$. In particular,
the constant function $1$ can be used as a barrier to show that
\[
||w||_{L^\infty (M)} \leq c \, ||{\mathbb L}_0 \,  w||_{L^\infty (M)} ,
\]
for some constant $c >0$ which does not depend on $\e$. Then, the
estimates for the derivatives of $w$ are consequences of Schauder's
estimates on geodesic balls of radius $\e$. \hfill $\Box$

Collecting the results of the previous sections we construct for the
operator ${\mathbb L}_{\e, N}$ a right inverse  whose norm is
uniformly bounded as $\e$ tends to $0$. The construction of the right
inverse relies on all our former analysis. To be more precise, we
glue together local parametrizes given by Proposition~\ref{pr:9.1}
and Lemma~\ref{le:fgfg} and obtain the right inverse for ${\mathbb
L}_{\e, N}$ by applying a perturbation argument together with
Proposition~\ref{co:11.2}. As will become clear in the proof and in
the statement of the next result, we need to assume that the
submanifold $N$ has one degree of regularity higher than would be
expected. Indeed, a natural guess would be that, in order to solve
the equation ${\mathbb L}_{\e, N} w = f \in {\cal C}^{k-2, \alpha}$
in some H\"older space ${\cal C}^{k, \alpha}$, one would need to
assume that the hypersurface $N$ is itself ${\cal C}^{k, \alpha}$.
However, our construction relies on the use of
Proposition~\ref{co:11.2} together with the use of the projection
operator $\Pi_{\e, N}$ which both require $N$ to be a ${\cal C}^{k+1,
\alpha}$ hypersurface. We now state the main technical result of our
paper~:
\begin{proposition}
Fix $\ell \geq 2$ and assume $N$ satisfies $(H)$ and $(H'_{\ell+1})$.
Then, there exists $\e_0 >0$ only depending on $\Lambda'$ and
$\Lambda''$ such that, for all $\e \in (0, \e_0)$, there exists an
operator
\[
{\mathbb G}_{\e, N} : {\cal E}^{\ell-2, \alpha}_{\e,N} (M)
\longrightarrow [{\cal E}^{\ell, \alpha}_{\e,N} (M)]_0,
\]
satisfying
\[
\Pi_{\e, N} \circ {\mathbb L}_{\e, N} \circ {\mathbb G}_{\e, N} = I .
\]
Furthermore the norm of ${\mathbb G}_{\e, N}$ is bounded by some
constant which depends on $\Lambda'$ and $\Lambda''$, but does not
depend on $\e \in (0, \e_0)$ nor on $N$.
\label{pr:th}
\end{proposition}
{\bf Proof~:} The construction of ${\mathbb G}_{\e, N}$ is decomposed
in $4$ steps. As already mentioned we "glue" together different
parametrizes which have been defined in the previous sections. In the
first step, we use the result of Lemma~\ref{le:fgfg} to reduce the
solvability of $\Pi_{\e, N} \circ {\mathbb L}_{\e, N} w =f$ to the
case where $f$ is supported in $V_{\tau_0}(N)$. This allows us to
use, in step 2 the result of Proposition~\ref{pr:9.1} where a right
inverse for the operator ${\mathfrak L}_\e$ has been constructed. In
step 3, we use the expansion of the Laplace-Beltrami operator in
Fermi coordinates which is provided by Proposition~\ref{co:11.2} to
estimate the difference between the operators ${\mathbb L}_{\e, N}$
and the operator ${\mathfrak L}_\e$.  At this point we have produced
a bounded operator $G$ which is almost a right inverse. In Step 4, it
will remain to apply a standard perturbation argument to find a right
inverse for $\Pi_{\e, N}\circ {\mathbb L}_{\e, N}$.

In the proof, the constant $c >0$ is a constant which may vary from
line to line, may depend on $\Lambda'$ and $\Lambda''$, but does not
depend on $\e$ (provided this parameter is chosen small enough), does
not depend on $f$, nor on $N$ chosen to satisfy $(H)$ and
$(H'_{\ell+1})$.

\noindent
{\bf Step 1.} Thanks to the result of Lemma~\ref{le:fgfg}, we find
$w_1 \in [{\cal C}^{\ell, \alpha}_\e (M)]_0$ solution of
${\mathbb L}_0 w_1 = f$. Furthermore, we know that
\[
||w_1||_{{\cal C}^{\ell, \alpha}_\e (M)} \leq c \, ||f||_{{\cal
C}^{\ell-2, \alpha}_\e (M)}.
\]
Now, the result of Lemma~\ref{le:proj} implies that we also have
\begin{equation}
||\Pi_{\e, N} \, w_1||_{{\cal C}^{\ell, \alpha}_\e (M)} \leq c \,
||f||_{{\cal C}^{\ell-2, \alpha}_\e (M)}.
\label{A}
\end{equation}
Observe that, as already mentioned, we already need to assume that
the hypersurface $N$ has bounded norm in ${\cal C}^{\ell+1,
\alpha}_\e (N_0)$.

We define the function
\[
g : =   f - {\mathbb L}_{\e, N} \circ \Pi_{\e, N} \, w_1 .
\]
Since, away from $V_{\tau_0}(N)$, we have $\Pi_{\e, N} u= u$ and
${\mathbb L}_0 =  {\mathbb L}_{\e, N}$, we conclude that, the
function $g$ is supported in $V_{\tau_0}(N)$. If we identify
$V_{\tau_0}(N)$ with $(-\tau_0, \tau_0 )\times N$ {\it via} (twisted)
Fermi coordinates, we can extend the function $g$ by $0$ to all
${\R}\times N$. We claim that
\begin{equation}
||g||_{{\cal C}^{\ell-2, \alpha}_{\gamma , \e}({\R}\times N)} \leq c
\, ||f||_{{\cal C}^{\ell-2, \alpha}_\e (M)} ,
\label{eq:resultat1}
\end{equation}
where $\gamma:= (\gamma_-, \gamma_+)$ are the indicial roots defined
in (\ref{eq:5.2}). Indeed, using (\ref{eq:5.3}) and (\ref{eq:5.4}),
we see that $u_\star$ converges exponentially to $\pm 1$ at a rate
which is dictated by the indicial roots $\gamma_\pm$. This in turn
shows that this is also the case for the difference between the
potential of the operators ${\mathbb L}_{\e, N}$ and ${\mathbb L}_0$.
Namely
\[
||\Gamma -  W''(u_\star(\cdot/\e))||_{{\cal C}^{\ell-2,
\alpha}_{\gamma, \e} ({\R}\times N)}\leq c .
\]
This inequality already implies that
\[
|| {\mathbb L}_{\e, N} \, w_1  - f||_{{\cal C}^{\ell-2,
\alpha}_{\gamma , \e}({\R}\times N)} \leq c \, ||f||_{{\cal
C}^{\ell-2, \alpha}_\e (M)}.
\]
Furthermore, using the definition of $\Pi_{\e, N}$ together with the fact that
\begin{equation}
|| \chi \, w_\star(\cdot/\e)||_{{\cal C}^{\ell, \alpha}_{\gamma, \e}
({\R}\times N)} \leq c,
\label{eq:estwstar}
\end{equation}
we get
\[
|| {\mathbb L}_{\e, N} \, (w_1 - \Pi_{\e, N} w_1)||_{{\cal
C}^{\ell-2, \alpha}_{\gamma , \e}({\R}\times N)} \leq c \,
||f||_{{\cal C}^{\ell-2, \alpha}_\e (M)}.
\]
This completes the proof of the claim.

\noindent
{\bf Step 2.} We define the projection operator
\[
\Pi^0_{\e, N} (u) (t,y) : = u (t,y) - \displaystyle
\frac{\displaystyle \int_{\R} u (s, y) \, w_\star (s /\e) \, ds
}{\displaystyle \int_{\R} w_\star^2 (s /\e) \, ds} \, w_\star (t/\e).
\]
Clearly,
\[
\Pi^0_{\e,N} :  {\cal C}^{\ell-2, \alpha}_{\gamma , \e} ({\R}\times
N)  \longrightarrow {\cal E}^{\ell-2, \alpha}_{\gamma , \e}
({\R}\times N) ,
\]
is bounded uniformly in $\e$. Hence, using (\ref{eq:resultat1}), we
conclude that
\begin{equation}
|| \Pi^0_{\e, N} \, g ||_{{\cal C}^{\ell-2, \alpha}_{\gamma ,
\e}({\R}\times N)} \leq c \, ||f||_{{\cal C}^{\ell-2, \alpha}_\e (M)}.
\label{eq:ששש}
\end{equation}
We use the result of Proposition~\ref{pr:9.1}, with some fixed
weights $\delta : = (\delta_-, \delta_+)$ satisfying  $- \gamma_\pm <
\delta_\pm <0$, to define $w_2 \in [{\cal E}^{\ell, \alpha}_{\delta,
\e} ({\R}\times N)]_{0}$ solution of
\[
{\mathfrak L}_\e \,  w_2 = \Pi^0_{\e, N} \, g ,
\]
where we recall that
\[
{\mathfrak L}_\e : =  -\e^2 (\del_t^2 + \Delta_N )+ \frac{1}{2} \,
W''(u_\star(\cdot/\e)) .
\]
Since the inverse of ${\mathfrak L}_\e$ has been shown to be
uniformly bounded, we get the estimate
\[
||w_2||_{{\cal C}^{\ell, \alpha}_{\delta, \e} ({\R}\times N)} \leq  c
\, ||\Pi_{\e, N}^0 \, g||_{{\cal C}^{\ell-2, \alpha}_{\delta, \e}
({\R}\times N)}\leq c \, ||\Pi_{\e, N}^0 \, g||_{{\cal C}^{\ell-2,
\alpha}_{\gamma, \e} ({\R}\times N)},
\]
which, together with (\ref{eq:ששש}) implies that
\begin{equation}
||w_2||_{{\cal C}^{\ell, \alpha}_{\delta, \e} ({\R}\times N)} \leq c
\, ||f||_{{\cal C}^{\ell-2, \alpha}_{\e} (M)}.
\label{eq:ppp}
\end{equation}

\noindent
{\bf Step 3.} We claim that
\begin{equation}
||\Pi_{\e, N} \, \left({\mathbb L}_{\e, N} \circ \Pi_{\e, N} \, (\chi
w_2) -  g \right) ||_{{\cal C}^{\ell-2, \alpha}_{\e} (M)} \leq \, c
\, \e \, || f||_{{\cal C}^{\ell-2, \alpha}_\e (M)} .
\label{wte}
\end{equation}

Since $w_2$ and $g$ are exponentially decaying in terms of $|t|/\e$,
we readily have
\[
||(1- \chi )  \, \Pi_{\e, N} \, \left({\mathbb L}_{\e, N} \circ
\Pi_{\e, N} \, (\chi w_2) - g \right)||_{{\cal C}^{\ell-2,
\alpha}_{\e} (M)} \leq \, c \, e^{-\kappa_1 / \e} \, || f||_{{\cal
C}^{\ell-2, \alpha}_\e (M)},
\]
for some fixed $\kappa_1 >0$ only depending on $\gamma$, $\delta$ and
$\tau_0$. Hence, identifying $V_{\tau_0}(N)$ with $[-\tau_0,
\tau_0]\times N$ as above, we only need to show that
\[
|| \chi  \, \Pi_{\e, N} \, \left({\mathbb L}_{\e, N} \circ \Pi_{\e,
N} \, (\chi w_2) -  g \right)||_{{\cal C}^{\ell-2, \alpha}_{\e}
({\R}\times N)} \leq \, c \, \e \, || f||_{{\cal C}^{\ell-2,
\alpha}_\e (M)} .
\]
The proof of this inequality essentially follows from the result of
Proposition~\ref{co:11.2} which gives the expansion of $\Delta_g$ in
Fermi coordinates relative to $N$. Observe that, again, we need to
assume that the hypersurface $N$ has bounded norm in ${\cal
C}^{\ell+1, \alpha}_\e(N_0)$ in order to apply the result of
Proposition~\ref{co:11.2} with the correct regularity.

To begin with we replace $\Pi_{\e, N}$ by $\Pi^0_{\e, N}$ in the left
hand side of (\ref{wte}). This introduce a discrepancy which we now
estimate. We set
\[
D_1 : =  \Pi_{\e, N} \circ {\mathbb L}_{\e, N} \circ \Pi_{\e, N} \,
(\chi w_2) - \Pi_{\e, N}^0 \circ {\mathbb L}_{\e, N} \circ \Pi_{\e,
N}^0 \, (\chi w_2) ,
\]
and
\[
D_2 : =  \Pi_{\e, N}^0 \, g -\Pi_{\e, N} \, g .
\]
Using the fact that $w_\star(\cdot/\e)$ is exponentially decaying in
terms of $|t|/\e$, we see that the same is true for $D_1$ and $D_2$,
so
\[
|| \chi \left(D_1 + D_2 \right)||_{{\cal C}^{\ell-2, \alpha}_{\e}
({\R}\times N)} \leq \, c \, e^{-\kappa_2 / \e} \, || f||_{{\cal
C}^{\ell-2, \alpha}_\e (M)} ,
\]
for some fixed positive $\kappa_2$ which only depends on $\gamma$
and $\tau_0$.

Having replaced $\Pi_{\e, N}$ by $\Pi_{\e,N}^0$, we now replace $\chi
\, w_2$ by $w_2$ in (\ref{wte}). This introduces yet another
discrepancy which we now estimate. We set
\[
D_3 : = \Pi_{\e, N}^0 \circ {\mathbb L}_{\e, N} \circ \Pi_{\e, N}^0
((1-\chi) w_2) .
\]
Since $w_2$ is exponentially decaying in terms of $|t|/\e$, we also have
\[
|| \chi \, D_3||_{{\cal C}^{\ell-2, \alpha}_{\e} ({\R}\times N)} \leq
\, c \, e^{-\kappa_3 /\e} \, || f||_{{\cal C}^{\ell-2, \alpha}_\e
(M)},
\]
for some fixed positive $\kappa_3$ which only depends on $\delta$ and $\tau_0$.

Having done the above modifications, we are left with the estimate of
\[
\begin{array}{rllll}
D_4  & : = & \Pi_{\e, N}^0 \circ {\mathbb L}_{\e, N} \circ \Pi_{\e,
N}^0 \, w_2 - \Pi_{\e, N}^0 \, g \\[3mm]
& = &  \Pi_{\e, N}^0 \circ {\mathbb L}_{\e, N} \,  w_2 - {\mathfrak
L}_\e  w_2 \\[3mm]
& = &  \Pi_{\e, N}^0 \left( {\mathbb L}_{\e, N} \,  w_2 - {\mathfrak
L}_\e  w_2 \right).
\end{array}
\]
We have used the fact that $\Pi^0_{\e, N} \, w_2 = w_2$ since $w_2
\in {\cal E}^{\ell, \alpha}_{\delta, \e} ({\R}\times N)$ and
${\mathfrak L}_\e \, w_2 =  \Pi^0_\e \, g$ to obtain the first
identity. In order to obtain the second identity, we have used the
fact that ${\mathfrak L}_\e  w_2 = \Pi_{\e, N}^0 \, {\mathfrak L}_\e
w_2 $, by construction of the inverse of ${\mathfrak L}_\e$. At this
point, we appeal for the result of Proposition~\ref{co:11.2} which
gives the expansion of $\Delta_g$ in Fermi coordinates about $N$ to
get
\[
{\mathbb L}_{\e, N} - {\cal L}_\e = L_1+ L_2 ,
\]
where the properties of the operators $L_1$ and $L_2$ are stated in
Proposition~\ref{co:11.2}. These properties, together with
(\ref{eq:ppp}), imply that
\[
|| \chi  \, D_4 ||_{{\cal C}^{\ell-2, \alpha}_{\e} ({\R}\times N)}
\leq \, c \, \e \, || f||_{{\cal C}^{\ell-2, \alpha}_\e (M)} .
\]
This completes the proof of the claim.

\noindent
{\bf Step 4.}  Collecting (\ref{wte}) together with the definition of
$g$, we conclude that
\[
|| \Pi_{\e, N} \circ {\mathbb L}_{\e, N} \circ \Pi_{\e, N}(w_1 + \chi
w_2) - f||_{{\cal C}^{\ell-2, \alpha}_{\e} (M)}\leq \, c \, \e \, ||
f||_{{\cal C}^{\ell-2, \alpha}_\e (M)} .
\]
Now, collecting (\ref{A}) and (\ref{eq:ppp}) together with
Lemma~\ref{le:proj}), we obtain
\[
|| \Pi_{\e, N} (w_1 + \chi w_2)||_{{\cal C}^{\ell, \alpha}_{\e} (M)}
\leq \, c \, || f||_{{\cal C}^{\ell-2, \alpha}_\e (M)}.
\]
At this point, the existence of ${\mathbb G}_{\e, N}$ easily follows
from a classical perturbation argument. \hfill $\Box$

Our next proposition gives an estimate of the norm of the operator
\[
(I-\Pi_{\e, N}) \circ {\mathbb L}_{\e, N} \circ {\mathbb G}_{\e, N} ,
\]
which, in some sense, measures the distance between the operator
${\mathbb L}_{\e, N} \circ {\mathbb G}_{\e, N}$ and the identity in
${\cal E}^{\ell-2, \alpha}_\e (M)$. We define the operator $S_{\e,
N}$ by
\begin{equation}
S_{\e, N} (u) (y) : =  \int_{\R} \chi  (s) \, u (s, y) \, w_\star (s
/\e) \, ds,
\label{eq:Sen}
\end{equation}
where $(t,y)$ are (twisted) Fermi coordinates relative to $N$ and $u$
is a function defined in $M$. Since $N$ is assumed to be a normal
graph over $N_0$, any function on $N$ can be identified with a
function on $N_0$ and its norm can be evaluated using the norms
defined in (\ref{eq:10.10}) or in (\ref{eq:10.1010}), provided $N$ is
regular enough. This being understood, we have~:
\begin{proposition}
Fix $\ell \geq 2$ and assume $N$ satisfies $(H)$ and $(H'_{\ell+1})$.
Then, there exists $c >0$ only depending on $\Lambda'$ and $
\Lambda''$ such that, for all $w \in [{\cal E}^{\ell, \alpha}_{\e,N}
(M)]_0$
\begin{equation}
| S_{\e, N} \circ {\mathbb L}_{\e, N}  \, w|_{{\cal C}^{\ell - 2,
\alpha}_\e (N_0)} \leq  c\, \e^2 \, ||w||_{{\cal C}^{\ell, \alpha}_\e
(M)} ,
\label{eq:ooo}
\end{equation}
where the norm $|\cdot|_{{\cal C}^{\ell', \alpha}_\e (N_0)}$ has been
defined in (\ref{eq:10.1010}).
\label{pr:ththth}
\end{proposition}
{\bf Proof~:} We have to estimate
\[
(t, y) \longrightarrow \int_{\R} \chi(s) \, w_\star (s/\e) L_{\e, N}
\, w (s,y) \, ds .
\]
We set
\[
E_1 : = \int_{\R} \chi(s) \, w_\star (s/\e) (L_{\e, N} - {\mathfrak
L}_\e) \, w (s,y) \, ds ,
\]
where ${\mathfrak L}_\e$ is defined as in Step 2 of the previous
proof. Using once more the result of Proposition~\ref{co:11.2}, we
obtain as in Step 3 of the previous proof
\[
|E_1 |_{{\cal C}^{\ell-2, \alpha}_\e (N_0)} \leq \, c \, \e^2 \, || w
||_{{\cal C}^{\ell, \alpha}_\e (M)} .
\]
We set
\[
E_2 : = \int_{\R}  w_\star (s/\e) \, [\chi , {\mathfrak L}_\e] \,
w(s,y)  \, ds ,
\]
where $[A, B]$ is the commutator of the operators $A$ and $B$. Using
the fact that $w_\star$ decays exponentially fast at $\pm \infty$, we
conclude that $E_2$ decays exponentially fast in terms of $|t|/\e$.
In particular
\[
|E_2|_{{\cal C}^{\ell-2, \alpha}_\e (N_0)} \leq \, c \,
e^{-\kappa_4/\e} \, || w ||_{{\cal C}^{\ell, \alpha}_\e (M)} ,
\]
for some fixed $\kappa_4 >0$ which only depends on $\gamma$ and $\tau_0$.

Having estimated $E_1$ and $E_2$, it remains to estimate
\[
y \longrightarrow  \int_{\R}  w_\star (s/\e) \, {\mathfrak L}_\e  (
\chi w) (s,y) \, ds .
\]
We set
\[
E_3 : = \int_{\R}  w_\star (s/\e) \, {\mathfrak L}_\e  \circ
(\Pi_{\e, N}-\Pi_{\e, N}^0) ( \chi w) (s,y) \, ds,
\]
and
\[
E_4 : = \int_{\R}  w_\star (s/\e) \, {\mathfrak L}_\e  \circ (I-
\Pi_{\e, N}) ( \chi w) (s,y) \, ds.
\]
Using the fact that $w_\star(\cdot/\e)$ is exponentially decaying in
terms of $|t|/\e$, we get
\[
|E_3 |_{{\cal C}^{\ell-2, \alpha}_\e (N_0)} \leq \, c \,
e^{-\kappa_5/\e} \, || u ||_{{\cal C}^{\ell, \alpha}_\e (M)},
\]
and, using in addition the fact that $(I-\Pi_{\e, N}) \, w = 0$, we also have
\[
| E_4 |_{{\cal C}^{\ell-2, \alpha}_\e (N_0)} \leq \, c \,
e^{-\kappa_5/\e} \, || u ||_{{\cal C}^{\ell, \alpha}_\e (M)},
\]
for some fixed $\kappa_5 >0$ which only depends on $\gamma$ and $\tau_0$.

Having estimated $E_1$ through $E_4$, it finally remains to estimate
\[
y \longrightarrow \int_{\R}  w_\star (s/\e) \, {\mathfrak L}_\e
\circ  \Pi_{\e, N}^0 ( \chi u) (s,y)  \, ds .
\]
But, as we have already done in (\ref{eq:XXX}), one checks that this
quantity is identically equal to $0$. This completes the proof of
(\ref{eq:ooo}). \hfill $\Box$

We would like to obtain similar results when $N$ is less regular than 
what is requires in the last two results. To do so we introduce 
smoothing operators which will alow us to find $N_\star$, smooth 
enough so that we can apply the previous results, and close enough to 
$N$ so that a perturbation argument can be applied.

We recall from \cite{Alin}, page 97, that there exists a one 
parameter family of smoothing operators $(R_\theta)_{\theta \geq 1}$ 
and $C >0$ such that
\begin{equation}
\begin{array}{rllll}
|| R_\theta \, u ||_{{\cal C}^{k, \alpha}} & \leq & C \, ||u 
||_{{\cal C}^{k', \alpha'}} & \quad \mbox{for} \quad k + \alpha \leq 
k'+ \alpha' \\[3mm]
|| R_\theta \, u ||_{{\cal C}^{k, \alpha}} & \leq & C \, \theta^{k+ 
\alpha - k'- \alpha'} \, ||u ||_{{\cal C}^{k', \alpha'}} & \quad 
\mbox{for} \quad k + \alpha \geq k'+ \alpha'\\[3mm]
|| u - R_\theta \, u ||_{{\cal C}^{k, \alpha}} & \leq & C \, 
\theta^{k+ \alpha - k'- \alpha'} \, ||u ||_{{\cal C}^{k', \alpha'}} & 
\quad \mbox{for} \quad k + \alpha \leq k'+ \alpha'.
\end{array}
\label{pr}
\end{equation}
These operators act on functions defined in ${\R}^n$ but they can be 
localized and extended to functions defined on smooth manifolds using 
a partition of unity.  We use $R_\theta$ to improve the regularity of 
a given hypersurface $N$ which is assumed to be a normal graph over 
$N_0$. We further assume that $(H)$ and $(H'_\ell)$ are satisfied. 
Hence, $N$ is the normal graph over $N_0$ for some function $\psi$ 
which satisfies
\[
||\psi||_{{\cal C}^{2}(N_0)} \leq \Lambda', \qquad \mbox{and}\qquad 
||\psi||_{{\cal C}^{\ell, \alpha}_\e (N_0)}\leq \Lambda'',
\]
for some $\ell \geq 2$. We define the function  $\psi_\star : = 
R_{1/\e}  \, \psi$, and we define the hypersurface $N_\star$ to be 
the graph of the function $\psi_\star$ over $N_0$. We claim that,
\begin{equation}
||\psi_\star ||_{{\cal C}^{2}(N_0)} \leq c \, ||\psi||_{{\cal 
C}^{2}(N_0)} , \qquad \qquad  |\psi_\star - \psi |_{{\cal C}^{\ell, 
\alpha}_\e (N_0)}\leq c \, \e^2 \, ||\psi||_{{\cal C}^{\ell, 
\alpha}_\e (N_0)},
\label{esd2}
\end{equation}
and, for all $\ell' \geq \ell$
\begin{equation}
||\psi_\star||_{{\cal C}^{\ell', \alpha}_\e (N_0)} \leq c \, 
||\psi||_{{\cal C}^{\ell, \alpha}_\e (N_0)},
\label{esd}
\end{equation}
where the constant $c$ depends on $\ell'$ and $\ell$ but does not 
depend on $\Lambda''$.

Indeed, using the first property of $R_\theta$ in (\ref{pr}), we 
readily obtain the first estimate of (\ref{esd2}). Similarly, we get 
$||\psi_\star||_{{\cal C}^{2} (N_0)}\leq c \, ||\psi||_{{\cal 
C}^{2}(N_0)}$ and, for all $2 \leq \ell' \leq \ell$,
\[
\begin{array}{rlllll}
||\psi_\star||_{{\cal C}^{\ell', \alpha} (N_0)} & \leq & c \, 
||\psi||_{{\cal C}^{\ell', \alpha}(N_0)}\\[3mm]
& \leq & c \, \e^{2- \ell'} \, ||\psi||_{{\cal C}^{\ell', \alpha}_\e(N_0)},
\end{array}
\]
which already proves that
\[
||\psi_\star||_{{\cal C}^{\ell, \alpha}_\e (N_0)} \leq c \, \Lambda''.
\]
Now, for all $\ell' \geq \ell$, we use the second property of 
$R_\theta$ in (\ref{pr}) to get
\[
\begin{array}{rlllll}
||\psi_\star||_{{\cal C}^{\ell', \alpha} (N_0)} & \leq & c \, 
\e^{\ell - \ell'}\, ||\psi||_{{\cal C}^{\ell, \alpha}(N_0)}\\[3mm]
& \leq & c \, \e^{2- \ell'} \, ||\psi||_{{\cal C}^{\ell, \alpha}_\e (N_0)},
\end{array}
\]
which implies the second estimate in (\ref{esd2}). Finally, for all 
$\ell' \leq \ell$, we use the third property of $R_\theta$ in 
(\ref{pr}) to show that
\[
\begin{array}{rlllll}
||\psi_\star - \psi||_{{\cal C}^{\ell', \alpha} (N_0)} & \leq & c \, 
\e^{\ell - \ell'}\, ||\psi||_{{\cal C}^{\ell, \alpha}(N_0)}\\[3mm]
& \leq & c \, \e^{2- \ell'} \, ||\psi||_{{\cal C}^{\ell, \alpha}_\e(N_0)},
\end{array}
\]
which proves (\ref{esd}). The proof of the claim is therefore complete.

In the case where $N$ has a nonempty boundary, this construction can 
be modified so that we can assume that $N_\star$ is an admissible 
hypersurface.  Thanks to the first estimate in (\ref{esd2}) and to 
(\ref{esd}), we can apply all the above results to $N_\star$. This 
yields the~:
\begin{proposition}
Fix $\ell \geq 2$ and $\tilde c >0$. Assume that $N$ satisfies $(H)$ 
and $(H'_\ell)$. Then, there exists $\e_0 >0$ only depending on 
$\tilde c$ such that, for all $2 \leq \ell' \leq \ell+2$ and for all 
$\e \in (0, \e_0)$, there exists an operator
\[
{\mathbb G}_{\e, N} : {\cal E}^{\ell'-2, \alpha}_{\e,N_\star } (M) 
\longrightarrow [{\cal E}^{\ell', \alpha}_{\e,N_\star} (M)]_0,
\]
satisfying
\[
\Pi_{\e, N_\star} \circ {\mathbb L}_{\e, N} \circ {\mathbb G}_{\e, N} = I .
\]
Furthermore,
\[
|| {\mathbb G}_{\e, N} \, w ||_{{\cal C}^{\ell', \alpha}_\e (M)} \leq 
c \,  ||w||_{{\cal C}^{\ell' - 2, \alpha}_\e (M)} ,
\]
and
\begin{equation}
| S_{\e, N_\star} \circ {\mathbb L}_{\e, N}  \circ {\mathbb G}_{\e, 
N} \, w |_{{\cal C}^{\ell' - 2, \alpha}_\e (N_0)} \leq  c\, \e^2 \, 
||w||_{{\cal C}^{\ell' - 2, \alpha}_\e (M)} ,
\label{kkc}
\end{equation}
for some constant $c>0$ which depends on $\Lambda'$ and $\tilde c$ 
but does not depend on $\e$ nor on $N$.
\label{pr:thth}
\end{proposition}
{\bf Proof~:}  We claim that
\begin{equation}
|| ({\mathbb L}_{\e, N_\star} - {\mathbb L}_{\e, N} ) \, w ||_{{\cal 
C}^{\ell'', \alpha}_\e (M)} \leq c \, \e \, ||w||_{{\cal C}^{\ell'', 
\alpha}_\e (M)} ,
\label{ssx}
\end{equation}
whenever $\ell''\leq \ell$.  This follows from the fact that the 
difference of these two operators is equal to the difference of their 
potentials which involve functions depending on the distance to $N$ 
and $N_\star$. The result is then a consequence of the fact that 
difference $\psi_\star -\psi$ satisfies the second estimate in 
(\ref{esd2}).

Now, we apply the results of Proposition~\ref{pr:th} and 
Proposition~\ref{pr:ththth} to get, for all $\ell' \geq 2$ and for 
all $f \in {\cal C}^{\ell'-2, \alpha}_\e (M)$ a solution $w \in 
[{\cal E}^{\ell', \alpha}_{\e, N_\star} (M)]_0$ of
\[
\Pi_{\e, N_\star} \circ {\mathbb L}_{\e, N_\star} \, w  = \Pi_{\e, N_\star} f ,
\]
with
\[
|| w ||_{{\cal C}^{\ell', \alpha}_\e (M)} \leq  c \, ||f||_{{\cal 
C}^{\ell' -2, \alpha}_\e (M)} .
\]
The existence of ${\mathbb G}_{\e, N}$ then follows from (\ref{ssx}) 
together with a simple perturbation argument. It remains to check 
(\ref{kkc}). But this follows at once from (\ref{ssx}) together with 
the result of Proposition~\ref{pr:ththth}. \hfill $\Box$

Observe that the mapping ${\mathbb G}_{\e, N}$ depends continuously 
on $N$. This essentially follows from the fact that, each step of the 
construction of this operator depends continuously on $N$.

\section{Moving the nodal set}\setcounter{equation}{0}

Assume that we are given an admissible hypersurface $N_0$. We define 
for convenience
\begin{equation}
{\mathbb Q}_\e (N) : = - \e^2 \, \Delta u_{\e, N} + \frac{1}{2} \, 
W'(u_{\e, N}),
\label{eq:pipo}
\end{equation}
where $N$ is an admissible hypersurface which can be written as a 
normal graph over $N_0$.  If $N_\star$ is the regularized 
hypersurface defined at the end of the previous section,  we define, 
as in (\ref{eq:Sen}), the operator $S_{\e, N_\star}$ by
\begin{equation}
S_{\e, N_\star} (u) (y) : =  \int_{\R} \chi  (s) \, u (s, y) \, 
w_\star (s /\e) \, ds,
\label{eq:13.9999}
\end{equation}
where $(s,y)$ are (twisted) Fermi coordinates relative to $N_\star$ 
and $u$ is a function defined in $V_{\tau_0}(N)$. Since $N_\star$ is 
assumed to be a normal graph over $N_0$, any function on $N_\star$ 
can be identified with a function on $N_0$ and its norm can be 
evaluated using the norms defined in (\ref{eq:10.10}) or 
(\ref{eq:10.1010}).

In this section we exploit the notion of nondegeneracy which have 
been defined in Definition~\ref{de:3.4} and Definition~\ref{de:3.44}. 
There are two results which correspond to the two different notions 
of non degeneracy. To start with let us assume that we are dealing 
with an admissible nondegenerate minimal hypersurface $N_0$. We have~:
\begin{lemma}
Assume that $\ell \geq 2$ is fixed. Assume that $N_0$ is an 
admissible nondegenerate minimal hypersurface. Then, there exists 
$\e_0 \in (0, 1)$ and $\bar c >0$ such that for all $f \in {\cal 
C}^{\ell-2, \alpha} (N_0)$ such that
\[
|f|_{{\cal C}^{\ell-2, \alpha}_\e (N_0)} \leq \bar c ,
\]
one can find an admissible hypersurface $N$ satisfying
\[
S_{\e, N_\star} \, {\mathbb Q}_\e (N) = \e^2 \, f
\]
Moreover, if $\psi_N$ is the function whose graph is $N$, we have
\[
||\psi_N||_{{\cal C}^{\ell, \alpha}_\e(N_0)} \leq c \, |f|_{{\cal 
C}^{\ell-2, \alpha}_\e (N_0)} ,
\]
for some constant $c>0$ which does not depend on $f$ nor on $\e$.
\label{le:min}
\end{lemma}
{\bf Proof~:} The proof of this result follows from the implicit 
function theorem. Granted the definition of $S_{\e, N_\star}$, we 
have to find and admissible hypersurface $N$ such that
\[
\ds \int_{\R} \chi (t) \, \left(   - \e^2 \, \Delta u_{\e, N}  + 
\frac{1}{2} \, W'(u_{\e,N}) \right) \, w_\star (t/\e) \, dt = \e^2 \, 
f ,
\]
where $(t,y)$ are Fermi coordinates relative to $N_\star$. Using 
(\ref{eq:11.2}) we compute
\[
- \e^2 \, \Delta u_{\e, N}  + \frac{1}{2} \, W'(u_{\e, N})  =   \e \, 
n \, H_{N_{s(t,y)}} \, w_\star (t/\e)
\]
where $s(t,y)$ denotes the distance to of the point $(t,y)$ to $N$ 
and $N_s$ is the hypersurface parallel to $N$ at distance $s$. We now 
write $N$ as a graph over $N_0$ for some function $\psi$ and hence 
$N_\star$ is the normal graph over $N_0$ for the function $\psi_\star 
= R_{1/\e}\, \psi$. We define
\[
A_\e (\psi) : =  n \, \ds \int_{\R} \chi (\e \, t) \, H_{N_{s(\e 
t,y)}} \, w_\star^2 (t) \, dt .
\]
And we claim that
\[
A_\e (\psi) = c_\star \, {\cal L}_{N_0} \, \psi + {\cal Q}(\e, \psi),
\]
where ${\cal L}_{N_0}$ is the Jacobi operator about $N_0$, where
\[
c_\star = \int_{\R} w_\star^2 \, dt,
\]
and where ${\cal Q}$ satisfies
\begin{equation}
| {\cal Q}(\e, \psi) |_{{\cal C}^{\ell-2, \alpha}_\e (N_0)} \leq c \, \e^2,
\label{za}
\end{equation}
and
\begin{equation}
| {\cal Q}(\e, \psi') - {\cal Q}(\e, \psi) |_{{\cal C}^{\ell-2, 
\alpha}_\e (N_0)} \leq c \, \e^2 \,  || \psi' - \psi ||_{{\cal 
C}^{\ell, \alpha}_\e (N_0)} ,
\label{zb}
\end{equation}
for some constant $c>0$ which depends on $\bar c$ but does not depend 
on $\e$ nor on $\psi, \psi'$ satisfying
\[
||\psi ||_{{\cal C}^{\ell, \alpha}_\e (N_0)} + ||\psi' ||_{{\cal 
C}^{\ell, \alpha}_\e (N_0)} \leq 1.
\]

To obtain this expansion, we first use the second estimate in 
(\ref{esd2}) to reduce to the case where $N_\star$ is replaced by 
$N$, i.e. $A_\e (\psi)$ is replaced by
\[
\tilde A_\e (\psi) : =  n \, \ds \int_{\R} \chi (\e \, s) \, H_{N_{\e 
s}} \, w_\star^2 (s) \, ds .
\]
If $D_\e (\psi) : = \tilde A_\e (\psi) - A_\e (\psi)$, we obtain from 
(\ref{esd2}) estimates for $D_\e$ which are similar to (\ref{za}) and 
(\ref{zb}).

We now use the expansion
\[
H_{N_t} =  H_{N} + t \, Q_1(\psi; y) + t^2 \, Q_2(\psi; t,y) ,
\]
where $Q_1(\psi; y)$ does not depend on $t$. In particular
\[
n \, \ds \int_{\R} \chi (\e \, s) \, \e \, s \, Q_1(\psi; y) \, 
w_\star^2 (s) \, ds \equiv 0,
\]
and, if we define
\[
E_{\e} (\psi): = n \, \e^2 \, \ds \int_{\R} \chi (\e \, s) \,  s^2 \, 
Q_2 (\psi; \e s, y) \, w_\star^2 (s) \, ds,
\]
we obtain for $E_\e$ estimates similar to (\ref{za}) and (\ref{zb}).

Finally, we use the fact that
\[
n \, H_N = n \, H_{N_0} + {\cal L}_{N_0} \psi + Q (\psi),
\]
and we check that
\[
F_{\e} (\psi): = n \, Q (\psi) \,  \ds \int_{\R} \chi (\e \, s) \, 
w^2_\star (s) \, ds,
\]
satisfies estimates similar to (\ref{za}) and (\ref{zb}).

We have assumed that the operator
\[
{\cal L}_{N_0} : [{\cal C}^{2, \alpha}(N_0)]_0 \longrightarrow {\cal 
C}^{0, \alpha} (N_0),
\]
is an isomorphism. Here the subscript $0$ means that we are 
considering functions which satisfy ${\cal B}_{N_0}(\psi)=0$ on $\del 
N_0$, if this boundary is not empty. This also implies that
\[
{\cal L}_{N_0} : [{\cal C}^{\ell, \alpha}_\e (N_0)]_0\longrightarrow 
{\cal C}^{\ell-2, \alpha}_\e (N_0),
\]
is an isomorphism, whose nor does not depend on $\e$ provided ${\cal 
C}^{\ell, \alpha}_\e (N_0)$ is endowed with the norm 
$||\cdot||_{{\cal C}^{\ell, \alpha}_\e (N_0)}$ defined in 
(\ref{eq:10.10}) and ${\cal C}^{\ell-2, \alpha}_\e (N_0)$ is endowed 
with the norm $|\cdot |_{{\cal C}^{\ell-2, \alpha}_\e (N_0)}$ defined 
in (\ref{eq:10.1010}). Indeed, if ${\cal L}_{N_0} w =f$ and if 
$|f|_{{\cal C}^{\ell-2, \alpha}_\e (N_0)} \leq 1$, then one can check 
that both $w$ and its derivatives up to second order are bounded 
functions on $N_0$ while the second partial derivatives of $w$ have 
their $|\cdot |_{{\cal C}^{\ell-2, \alpha}_\e (N_0)}$ norm bounded.

The proof of the result is now a corollary of the fixed point theorem 
for contraction mapping, provided the constant $\bar c$ is chosen 
small enough. \hfill $\Box$

We now state the corresponding result for volume-nondegenerate 
constant mean curvature hypersurfaces.
\begin{lemma}
Assume that $\ell \geq 2$is fixed. Assume that $N_0$ is an admissible 
constant mean curvature hypersurface which is volume-nondegenerate. 
We define
\[
\lambda_\star : = \frac{1}{2}  \, c_\star \, n \, H_{N_0}
\]
where
\[
c_\star : = \int_{-1}^1 \sqrt{W(s)} \, ds .
\]
Then, there exists $\e_0 \in (0, 1)$ and $\bar c >0$ such that for 
all $f \in {\cal C}^{\ell, \alpha} (N_0)$ and all $\mu \in {\R}$, 
satisfying
\[
|f|_{{\cal C}^{\ell-2, \alpha}_\e (N_0)} \leq \bar c , \qquad 
\mbox{and} \qquad |\mu|\leq \bar c ,
\]
one can find an admissible hypersurface $N$ and a constant $\lambda$ such that
\[
S_{\e, N_\star} \,  \left({\mathbb Q}_\e (N) - \e \, \lambda \right) 
=  \e^2\, f ,
\]
and
\[
\int_M u_{\e, N} \, dv_g = c_0 \, |M| + \mu.
\]
Moreover, if $\psi_N$ is the function whose graph is $N$, we have
\[
||\psi_N||_{{\cal C}^{\ell, \alpha}_\e(N_0)} \leq c \, (|f|_{{\cal 
C}^{\ell-2, \alpha}_\e (N_0)} + |\mu|),
\]
for some constant $c>0$ which does not depend on the data $f$,  $\mu$ 
nor on $\e$.
\label{le:cmc}
\end{lemma}
{\bf Proof~:} We now  define
\[
A_\e (\psi) : =  n \, \ds \int_{\R} \chi (\e \, t) \, H_{N_{s(\e 
t,y)}} \, w_\star^2 (t) \, dt - \lambda \, \int_{\R}\chi (\e \, t) \, 
w_\star(t) \, dt ,
\]
which can be expanded as
\[
A_\e (\psi) = \ds  n \, H_{N_0} \, \int_{\R}  w^2_\star (t) \, dt - 
\e \, \lambda \, \int_{\R}  w_\star(t) \, dt +  c_\star \, {\cal 
L}_{N_0} \, \psi + {\cal Q}(\e, \psi) ,
\]
where ${\cal Q}(\e, \psi)$ satisfies (\ref{za}) and (\ref{zb}). Now, 
observe that
\[
\int_{\R} w_\star \, dt= 2 ,
\]
and also that
\[
\int_{\R} w^2_\star \, dt =  \int_{\R}\sqrt{W'(u_\star)} \, \del_t 
u_\star \, dt = \int_{-1}^1  \sqrt{W'(s)} \, ds = c_\star.
\]
Granted the definition of $\lambda_\star$, we can write
\[
A_\e (\psi) =  c_\star \, {\cal L}_{N_0} \, \psi  - 2 \, (\lambda- 
\lambda_\star) + {\cal Q}(\e, \psi) ,
\]

We now expand
\[
\int_M u_{\e, N} \, dv_g =  c_0 \, |M| - 2 \int_{N_0} \, \psi \, da_g 
+ O(\e, \psi),
\]
where the operator $O$ satisfies
\[
|O(\e, \psi) |\leq c \, \e^2,
\]
and
\[
| O(\e, \psi') - O(\e, \psi) |  \leq c \, \e^2 \,  || \psi' - \psi 
||_{{\cal C}^{\ell, \alpha}_\e (N_0)} ,
\]
for some constant $c>0$ which depends on $\bar c$ but does not depend 
on $\e$ nor on $\psi, \psi'$ satisfying
\[
||\psi ||_{{\cal C}^{\ell, \alpha}_\e (N_0)} + ||\psi' ||_{{\cal 
C}^{\ell, \alpha}_\e (N_0)} \leq 1.
\]

The proof now follows exactly the proof of the previous Lemma when 
volume-nondegeneracy replaces nondegeneracy. We omit the details. 
\hfill $\Box$

\section{The nonlinear problem}\setcounter{equation}{0}

\subsection{The proof of Theorem 1}

Assume that we are given a volume-nondegenera\-te admissible minimal 
hypersurface $N_0$ in $M$. We would like to solve the nonlinear 
problem
\begin{equation}
- \e^2 \, \Delta_g (u_{\e, N} + v) + \frac{1}{2} \, W'(u_{\e, N} + v ) =0 ,
\label{eq:13.31}
\end{equation}
in $M$. In addition, we want $u_{\e, N}+v$ to have $0$ Neumann 
boundary data on $\del M$ if this later is not empty. This means 
that, for all $\e$ small enough, we would like to find an admissible 
hypersurface $N$ close to $N_0$ and a function $v$ close to $0$ 
satisfying (\ref{eq:13.31}).

Recall that we have defined the nonlinear operator
\begin{equation}
{\mathbb Q}_\e (N) : = - \e^2 \, \Delta_g u_{\e, N}  + \frac{1}{2} \, 
W'(u_{\e, N}),
\label{eq:jkl}
\end{equation}
which corresponds to the error we produce when we consider $u_{\e, 
N}$ as a solution of (\ref{eq:13.31}).  We define
\begin{equation}
\tilde {\mathbb Q}_{\e} (N, v) : = \frac{1}{2} \, \left( W'(u_{\e, 
N}+v) - W'(u_{\e, N})  - W''(u_{\e, N}) \, v\right),
\label{eq:jkll}
\end{equation}
which is nothing but the Taylor expansion of the nonlinearity $W'$ at 
$u_{\e, N}$.

These definitions being understood, the equation we have to solve reads
\[
{\mathbb Q}_\e (N)  + {\mathbb L}_{\e, N} \, v + \tilde {\mathbb 
Q}_\e (N, v)  = 0 .
\]
Using the projection $\Pi_{\e, N_\star}$ defined in (\ref{eq:13.999}) 
and the operator $S_{\e, N_\star}$ defined in (\ref{eq:13.9999}), we 
conclude that (\ref{eq:13.31}) is equivalent to the system
\begin{equation}
\left\{
\begin{array}{rllll}
\ds \Pi_{\e, N_\star} \circ {\mathbb L}_{\e, N} \, v  & = & - 
\Pi_{\e, N_\star} \, \left( {\mathbb Q}_\e (N) + \tilde {\mathbb 
Q}_\e (N, v)  \right) \\[3mm]
\ds S_{\e, N_\star} \, {\mathbb Q}_\e (N) & = & - S_{\e, N_\star} \, 
\left( {\mathbb L}_{\e, N} \, v + \tilde {\mathbb Q}_{\e} (N,  v) 
\right) ,
\end{array}
\right.
\end{equation}
where $N_\star$ is the hypersurface obtained at the end of \S  11.

This system will be solved using a fixed point argument. The key 
result which allows one to apply a fixed point theorem is the 
following estimate~:
\begin{lemma}
There exists a constant $c_0 >0$ such that, for all $\e \in (0,1)$
\begin{equation}
|| {\mathbb Q}_\e (N_0)||_{{\cal C}^{0, \alpha}_\e(M)} \leq \ds 
\frac{c_0}{2} \, \e^2,
\label{eq:c_0}
\end{equation}
and, given $c_2 >0$, there exists $\e_0 >0$ such that, for all $\e 
\in (0,\e_0)$ and all admissible hypersurface $N$ satisfying
\[
||N||_{{\cal C}^{2, \alpha}_\e (N_0)} \leq c_2 \, \e^2 ,
\]
we have
\begin{equation}
||{\mathbb Q}_\e (N)||_{{\cal C}^{0, \alpha}_\e (M)} \leq c_0 \, \e^2.
\label{eq:c_00}
\end{equation}
\label{le:c_0}
\end{lemma}
{\bf Proof~:} It follows from (\ref{eq:11.2}) that
\[
{\mathbb Q}_\e (N) : = - \e^2 \, \Delta_g u_{\e, N}  + \frac{1}{2} \, 
W'(u_{\e, N}) = n \,  \e \, H_{N_{t}} \, w_\star (t/\e),
\]
in $V_{\tau_0/2}(N)$, where $N_{t}$ is the hypersurface parallel to 
$N$ at distance $t$. Moreover ${\mathbb Q}_\e (N) \equiv 0$ in 
$M-V_{\tau_0}(N)$. Since $N_0$ is minimal and $||N||_{{\cal C}^{2, 
\alpha}_\e (N_0)} \leq c_2 \, \e^2$, we estimate
\[
|H_{N_{t}}| \leq c \, ( |t|  + c_2 \, \e^2),
\]
$V_{\tau_0} (N)-V_{\tau_0/2}(N)$ and this already implies that
\[
||H_{N_{t}} \, w_\star (t/\e)||_{{L^\infty}(M)}\leq c \, (\e + c_2 \, \e^2).
\]
(The estimate in $V_{\tau_0} (N)-V_{\tau_0/2}(N)$ is easy to get 
since  ${\mathbb Q}_\e (N) $ is exponentially small in this set). 
The estimates for the H\"older derivative follows similarly.

It now suffices to apply this estimate to $N_0$ itself to obtain 
(\ref{eq:c_0}), while (\ref{eq:c_00}) follows at once by taking $\e$ 
to be small enough. \hfill $\Box$

Assume that we are given $v \in {\cal C}^{0, \alpha}_\e (M)$ and $N$ 
an admissible hypersurface close to $N_0$. To make things 
quantitatively precise, we assume that
\begin{equation}
||v||_{{\cal C}^{0, \alpha}_\e (M)} \leq c_1 \, \e^2 ,
\label{eq:c_1}
\end{equation}
and
\begin{equation}
||N||_{{\cal C}^{2, \alpha}_\e (M)} \leq c_2 \, \e^2 .
\label{eq:c_2}
\end{equation}
For some constant $c_1$ and $c_2$ which will be fixed shortly. We 
apply the result of Proposition~\ref{pr:thth} and Lemma~\ref{le:min} 
to find $\tilde v$ and $\tilde N$ solutions of
\begin{equation}
\left\{
\begin{array}{rllll}
\ds \Pi_{\e, N_\star} \circ {\mathbb L}_{\e, N} \, \tilde v  & = & - 
\Pi_{\e, N_\star} \, \left( {\mathbb Q}_\e (N) + \tilde {\mathbb 
Q}_\e (N, v)  \right) \\[3mm]
\ds S_{\e, \tilde N_\star} \, {\mathbb Q}_\e (\tilde N) & = & - 
S_{\e, N_\star} \, \left( {\mathbb L}_{\e, N} \, v + \tilde {\mathbb 
Q}_{\e} (N,  v) \right) .
\end{array}
\right.
\end{equation}
Using the fact that $\tilde Q_\e$ is quadratic in $v$ together with 
Lemma~\ref{le:c_0}, we easily get the estimates
\[
||\tilde v||_{{\cal C}^{2, \alpha}_\e (M)} \leq  c \, (c_0 \, \e^2 + 
c_1^2 \, \e^4 ),
\]
and
\[
|| \tilde N||_{{\cal C}^{2, \alpha}_\e (M)} \leq c \, (c_1 \, \e^2 + 
c_1^2 \, \e^4) .
\]
For some constant $c >0$ which does not depend on $c_1$ nor on $c_2$.

It follows at once that, if
\[
c_1 = 2 \, c\, c_0 \qquad \mbox{and} \qquad c_2 = 2 \, c \, c_1,
\]
and if $\e$ is chosen small enough, this produces a continuous 
mapping from the set of $(v, N) \in {\cal C}^{2, \alpha}_\e (M) 
\times {\cal C}^{2, \alpha}_\e (N_0)$ satisfying (\ref{eq:c_1}) and 
(\ref{eq:c_2}) into itself. If this mapping were compact, we would 
obtain a fixed point through Schauder's fixed point theorem. However, 
in our case the mapping just fails to be compact since no regularity 
is gained through the iteration process. To overcome this problem, we 
use once more the smoothing operators $R_\theta$ which have been 
introduced in the proof of Proposition~\ref{pr:thth}. We define 
instead $\tilde v$ and $\tilde N$ to be the solutions of
\begin{equation}
\left\{
\begin{array}{rllll}
\ds \Pi_{\e, N_\star} \circ {\mathbb L}_{\e, N} \tilde v  & = & - 
\Pi_{\e, N_\star} \circ R_\theta \,  \left( {\mathbb Q}_\e (N) + 
\tilde {\mathbb Q}_\e (N, v)  \right) \\[3mm]
\ds S_{\e, \tilde N_\star} \, {\mathbb Q}_\e (\tilde N) & = & - 
R_\theta \circ S_{\e, N_\star} \, \left( {\mathbb L}_{\e, N} \, v + 
\tilde {\mathbb Q}_{\e} (N,  v) \right) .
\end{array}
\right.
\label{zz}
\end{equation}
We fix $\alpha' >\alpha$ and choose $\theta$ so that $C \, 
\theta^{\alpha'- \alpha} = 2$, where $C$ is the constant which 
appears in the estimate of the smoothing operator $R_\theta$. If $c_1 
= 4 \, c \, c_0$  and $c_2 = 4 \, c \, c_1$, this produces, for all 
$\e$ small enough, a continuous mapping from the set of $(v, N) \in 
{\cal C}^{2, \alpha}_\e (M) \times {\cal C}^{2, \alpha}_\e (N_0)$ 
satisfying (\ref{eq:c_1}) and (\ref{eq:c_2}) into itself, but this 
tile the mapping is compact. We conclude that there exists a fixed 
point $(v_\theta$, $N_\theta)$.  Finally, we pass to the limit as 
$\theta$ tends to $+\infty$ (i.e. as $\alpha'$ tends to $\alpha$). 
The solutions $(v_\theta$, $N_\theta)$ of (\ref{zz}) being uniformly 
bounded in ${\cal C}^{2, \alpha}_\e (M) \times {\cal C}^{2, 
\alpha}_\e (N_0)$, we can extract a subsequence which converges to 
$(v, N)$ in ${\cal C}^{2, \beta}_\e (M) \times {\cal C}^{2, \beta}_\e 
(N_0)$, for some fixed $\beta < \alpha$. The limit $(v,N)$ is then a 
solution of our problem.

\subsection{The proof of Theorem 2}

Assume that we are given a volume-nondegenera\-te admissible constant 
mean curvature hypersurface $N_0$ in $M$. We would like to solve the 
nonlinear problem
\begin{equation}
- \e^2 \, \Delta (u_{\e, N} + v) + \frac{1}{2} \, W'(u_{\e, N} + v ) 
= \e \, \lambda ,
\label{eq:13.3}
\end{equation}
in $M$, with $u_{\e, N}+v$ having $0$ Neumann boundary data if $\del 
M$ is not empty. This equation has to be complimented with the 
constraint
\begin{equation}
\int_{M} (u_{\e, N} + v) \, dv_g = c_0 \, |M|,
\label{eq:13.5}
\end{equation}
where the constant $c_0 \in (-1,1)$ is fixed so that
\[
c_0 \, |M| =  |M^+(N_0)|-|M^-(N_0)|.
\]

Again, in order to solve (\ref{eq:13.3}), we write the equation as a 
fixed point problem
\begin{equation}
\left\{
\begin{array}{rllll}
\ds \Pi_{\e, N_\star} \circ {\mathbb L}_{\e, N} \, v  & = & \ds - 
\Pi_{\e, N_\star} \, \left( {\mathbb Q}_\e (N) - \e \, \lambda + 
{\mathbb Q}_\e (N, v)  \right) \\[3mm]
\ds S_{\e, N_\star} \, \left( {\mathbb Q}_\e (N) - \e \, \lambda 
\right) & = & - S_{\e, N_\star} \, \left( {\mathbb L}_{\e, N} \, v + 
\tilde {\mathbb Q}_{\e} (N,  v) \right) \\[3mm]
\ds \int_{M} u_{\e, N}\, dv_g & = & c_0 \, |M| - \ds \int_{M} \, v\, dv_g .
\end{array}
\right.
\end{equation}
The proof of the existence of a fixed point is identical to the proof 
of the previous result, with Lemma~\ref{le:min} replaced by 
Lemma~\ref{le:cmc}. The only difference being that the hypersurface 
$N_0$ does not have $0$ mean curvature anymore and this implies that 
Lemma~\ref{le:c_0} has to be replaced by
\begin{lemma}
There exists a constant $c_0 >0$ such that, for all $\e \in (0,1)$
\[
|| {\mathbb Q}_\e (N_0)||_{{\cal C}^{0, \alpha}_\e(M)} \leq \ds 
\frac{c_0}{2} \, \e,
\]
and, given $c_2 >0$, there exists $\e_0 >0$ such that, for all $\e 
\in (0,\e_0)$ and all admissible hypersurface $N$ satisfying
\[
||N||_{{\cal C}^{2, \alpha}_\e (N_0)} \leq c_2 \, \e ,
\]
we have
\[
||{\mathbb Q}_\e (N)||_{{\cal C}^{0, \alpha}_\e (M)} \leq  c_0 \, \e.
\]
\end{lemma}

The fact that we do not get an estimate as good as the one obtained 
in Lemma \ref{le:c_0} is a consequence of the fact that the mean 
curvature of $N_0$ is not necessarily equal to $0$.

This implies that (\ref{eq:c_1}) and (\ref{eq:c_2}) have to be 
replaced respectively by
\[
||v||_{{\cal C}^{0, \alpha}_\e (M)} \leq c_1 \, \e ,
\]
and
\[
||N||_{{\cal C}^{2, \alpha}_\e (M)} \leq c_2 \, \e .
\]
Details are left to the reader.

\providecommand{\bysame}{\leavevmode\hbox to3em{\hrulefill}\thinspace}

\end{document}